\documentclass{conm-p-l}
 
\usepackage{amsmath}
\usepackage{amsthm}
\usepackage{amssymb}
 
\newtheorem{thm}{Theorem}[section]
\newtheorem{cor}[thm]{Corollary}
\newtheorem{lem}[thm]{Lemma}
\newtheorem{prop}[thm]{Proposition}

\newtheorem{defn}[thm]{Definition}
\newtheorem{example}[thm]{Example}

 \newtheorem{notation}[thm]{Notation}

\newtheorem{rem}[thm]{Remark}

\theoremstyle{remark}

\numberwithin{equation}{section}

\newcommand{\R}{{\mathbb R}}
\newcommand{\N}{{\mathbb N}}
\newcommand{\Z}{{\mathbb Z}}
\newcommand{\C}{{\mathbb C}}
\newcommand{\T}{{\mathbb T}}

\newcommand{\rn}{{\mathbb R}^n}

\newcommand{\rp}{{\mathbb R}_+}

\newcommand{\A}{{\mathcal A}}
\newcommand{\B}{{\mathcal B}}

\newcommand{\cL}{{\mathcal L}}
\newcommand{\cC}{{\mathcal C}}
\newcommand{\D}{{\mathcal D}}

\newcommand{\cP}{{\mathcal P}}

\newcommand{\cS}{{\mathcal S}}

\newcommand{\cH}{{\mathcal H}}

\newcommand{\ga}{\alpha}
\newcommand{\gb}{\beta}
\renewcommand{\gg}{\gamma}
\newcommand{\gG}{\Gamma}

\newcommand{\gve}{\varepsilon}

\newcommand{\gl}{\lambda}
\newcommand{\gL}{\Lambda}
\newcommand{\go}{\omega}
\newcommand{\gO}{\Omega}

\newcommand{\gvp}{\varphi}
\newcommand{\gt}{\theta}

\newcommand{\gs}{\sigma}

\newcommand{\Op}{{\rm op\,}}

\newcommand{\n}[1]{\left<#1\right>}

\newcommand{\Proof}{\noindent{\em Proof}}

\newcommand{\eproof}{{~\hfill$ \triangleleft$}}
\newcommand{\ti}{\tilde}
\renewcommand{\i}{\infty}

\def\bli{\begin{enumerate}}
\def\eli{\end{enumerate}}
\def\vec#1#2{\begin{array}{c}#1\\#2\end{array}}
\def\vectd#1#2#3{\begin{array}{c}#1\\#2\\#3\end{array}}
\begin{document}
 
\title[Noncommutative Residues and  Trace  Expansions]{
Noncommutative Residues, Dixmier's Trace, and Heat Trace
Expansions on Manifolds with Boundary}
% author one information
\author{Elmar Schrohe}
\address{Institut f\"ur Mathematik, Universit\"at Potsdam, 14415 Potsdam, Germany}
\curraddr{}
\email{schrohe@math.uni-potsdam.de}
\thanks{This paper appeared in:
    B. Booss-Bavnbek, K. Wojciechowski (eds)
    Geometric Aspects of Partial Differential Equations.
    Proceedings of a Minisymposium on Spectral Invariants,
    Heat Equation Approach, September 18-19, 1998, Roskilde, Denmark.
    Contemporary Mathematics Vol 242,  pp. 161--186,
    Amer. Math. Soc., Providence, R.I., 1999.
}

%    General info
\date{}

\begin{abstract}
For manifolds with boundary, we define an extension of Wodzicki's noncommutative
residue to boundary value problems in Boutet de Monvel's calculus.
We show that this residue can be recovered with the help of heat kernel expansions 
and explore its relation to Dixmier's trace.
\end{abstract}

\maketitle
\tableofcontents
\section*{Introduction}

In 1984, M. Wodzicki discovered a trace on the algebra $\Psi_{cl}(\gO)$ of all classical pseudodifferential operators on a closed compact manifold $\gO$, \cite{WodThesis}. He called it the 
{\em noncommmutative residue}; meanwhile the notion {\em Wodzicki residue} is 
also widely used. This trace vanishes if the order of the operator is less than $-{\rm dim}(\gO)$.
For one thing this shows that the noncommutative residue is not an extension of the 
usual operator trace; it also implies that it is zero on the ideal $\Psi^{-\infty}(\gO)$
of all regularizing operators and therefore yields a trace on $\Psi_{cl}(\gO)/\Psi^{-\infty}(\gO)$. In fact it turns out to be the unique trace on this 
algebra up to multiples. 

Wodzicki's residue has found a wide range of applications both in mathematics
and mathematical physics. It plays a prominent role, for example, in Connes' noncommutative geometry. This is mainly due to the fact that it coincides with Dixmier's
trace on pseudodifferential operators of order $-{\rm dim}(\gO)$ as observed by Connes
\cite{Connes88}. Certain abstract expressions of noncommutative geometry thus can 
be computed explicitly for the classical differential geometric situation.
Moreover, while Dixmier's trace in general depends on the choice of an averaging 
procedure, this is not the case for these operators. 

The noncommutative residue is closely related to zeta functions of operators and
generalized heat trace asymptotics.
In fact this is how Wodzicki originally defined it:
Given a pseudodifferential operator $P$,
one may choose an invertible pseudodifferential
operator $A$ of order larger than that of $P$ and consider, for small $u$,
the zeta function $\zeta_{A+uP}(s)$.  
According to  classical results of Seeley, this is a holomorphic function
for large ${\rm Re}\,s$;
it extends meromorphically to the whole complex plane.
The formula
${\rm res}\,P = {\rm ord}\,A \frac d{du}
\left({\rm Res}_{s=1}\zeta_{A+uP}(s)\right)|_{u=0}$
then defines Wodzicki's residue (note that the right hand side is independent
of the choice of the auxiliary operator $A$).

Standard formulas link the zeta function of an operator $A$
and the so-called `heat trace' $\,{\rm trace}\,e^{-tA}$. The analysis
of both objects is therefore similar, and one obtains corresponding
results for the noncommutative residue in terms of coefficients
of heat trace expansions.
For practical purposes, however, an expression in terms of the symbol of $P$,
also derived by Wodzicki, is often more convenient, cf.~Theorem \ref{251.3}.

The present survey focuses on manifolds with boundary.
Wodzicki noticed right away that there is no trace on the algebra of
classical pseudodifferential symbols whenever the underlying manifold
is noncompact or has a boundary.
One might argue, however, that this algebra is not the natural object
to consider, since the standard symbol composition does not `feel' the presence of the
boundary.

The natural analogue of the algebra of pseudodifferential operators
on a closed
manifold is rather Boutet de Monvel's algebra of boundary value problems.
And
indeed, it turned out that on this algebra one can find a trace which is an
extension of Wodzicki's in the sense that both coincide when the boundary is 
empty, cf.~Fedosov, Golse, Leichtnam, Schrohe \cite{FGLS2}. It is therefore
also called a noncommutative residue. It vanishes on operators of sufficiently
low order, thus carries over to the quotient modulo regularizing elements. Very
much like in the classical case it turns out to be the unique trace;
for technical reasons, continuity is required here.

It is a natural question whether this trace is related to Dixmier's.
In fact, 
one can show that operators in Boutet de Monvel's calculus of order
$-\dim \Omega$  also belong to Dixmier's ideal.
It is a little surprising that both traces do not coincide.
It is possible, however, to prove a formula that computes Dixmier's
trace in terms of the symbols of the operators;
the ingredients in this formula
are the same as in the case of the noncommutative residue.
In particular, the
expression is local, and Dixmier's trace is independent of the averaging
procedure, cf.~Nest and Schrohe \cite{NestSchrohe98}.

It has been an open problem to link the noncommutative residue for
manifolds with boundary to heat trace expansions. Of course,
asymptotics for ${\rm trace}\,e^{-tB}$, where $B$ is an elliptic boundary
value problem in Boutet de Monvel's calculus, had been established
for example by Grubb in her 1986 book \cite{Grubb86}. Due to
technical difficulties, however, no complete expansion could be obtained.
The remainder term was $O(t^{1-\gve})$ for some $\gve>0$, and the
noncommutative residue therefore could not be detected.

In connection with their work on Atiyah-Patodi-Singer boundary
problems, Grubb and Seeley \cite{GrubbSeeley93a} overcame similar difficulties
by introducing a parameter-dependent pseudodifferential calculus
based on so-called weakly parametric symbols.
It recently became clear that this technique could also be applied here.
In joint work with Grubb, we fix a class of `nice' auxiliary
second order boundary
value problems, say $A$, and derive an asymptotic expansion
for ${\rm trace}\, Pe^{-tA}$ for every operator $P$ in Boutet de Monvel's
calculus. Moreover, we show how the coefficient of the
logarithmic term in this expansion relates to the noncommutative residue.

The paper is organized as follows. For the benefit of the reader who
is not familiar with the subject we start with a survey of the
situation in the boundaryless case.
There follows a short introduction to Boutet de Monvel's calculus
and the noncommutative residue on manifolds with boundary.
Next we analyze the relation to Dixmier's trace and give the
trace formula.
The subsequent section explains how the weakly parametric calculus
can be used to obtain heat trace asymptotics and the link between
the noncommutative residue and the logarithmic term. The final part gives
references to various applications and related work. 

\section{The Classical Results for Closed Manifolds}
This section contains a short introduction to the theorems of
Wodzicki and Connes for the boundaryless case.
In addition, the calculus of Grubb and Seeley \cite{GrubbSeeley93a}
and its relation to Wodzicki's residue are sketched.
Let us first fix the notation.

\begin{defn}\label{250.1}{\rm Let 
${\A}$ be an algebra over $\C$. A linear map $\tau: \A \to \C$  is called a {\it trace},
if it vanishes on commutators, i.e., if
$$\tau [P,Q] =
\tau (PQ-QP) =0 \mbox{~~for all~~} P,Q \in \A.$$
Clearly, if $\tau$ is a trace, then $\gl\tau$ is a trace for each $\gl$ in $\C$;
moreover, the zero map is always a trace. When we speak of a unique trace,
we shall mean that it is non-zero and the only one up to multiples.   
}\end{defn}

\begin{example}\label{250.2}{\rm On $M_r (\C)$, the algebra of $r\times r$ matrices 
over $\C$, there is a unique trace, namely the standard one, 
${\rm Tr}: A \mapsto \sum^r_j A_{jj}$.
Indeed, let  $E_{jk}$ denote the matrix
having a single 1 at position  $j,k$ (and zeros else). Then the statement is immediate 
from the observation that $[E_{jk},E_{kk}] = E_{jk}$ for $j\not= k$ and 
$[E_{jk},E_{kj}]= E_{jj} - E_{kk} $. 
}\end{example}

\notation\label{250.4}{\rm  $\gO$ is an $n$-dimensional compact manifold,
$E$ a vector bundle,
and  $\Psi_{cl}(\gO)$  the algebra of all classical pseudodifferential
operators on $\gO$,  acting on sections of $E$.
By   
$\Psi^{-\i}(\gO)$ we denote the ideal of regularizing elements
and by $\A$ the quotient $\Psi_{cl}(\gO)/\Psi^{-\i}(\gO) $.

Let $A\in \Psi_{cl}(\gO)$ be of order $m$.
Over each coordinate neighborhood $U$ its symbol $a$
has an asymptotic expansion $a\sim\sum_{j=0}^\i a_{m-j}$
into terms $a_{m-j}
= a_{m-j}(x,\xi)\in S^{m-j}(\R^n\times\R^n)$
that are homogeneous in $\xi$ of degree $m-j$ for $|\xi|\ge 1$.
Changing these terms smoothly for small $\xi$ results in regularizing terms.
The equivalence class of $A$ in  $\A$ therefore
can be identified with a formal sum of homogeneous functions
$
\sum^\i_{j=0} a_{m-j} (x,\xi),
$
where now the  $a_{m-j} \in C^\i (U \times (\R^n \setminus \{ 0 \}))$
are homogeneous in $\xi$ of degree $m-j$ taking values in square matrices. 
This is what we shall do in the following. }
\bigskip

\refstepcounter{thm}
{\sc \thesection.\arabic{thm}} \label{251.2}{\sc The Form   $\gs$.}
{\rm On $\R^n,~ n \geq 2$, 
define the $(n-1)$-form
$$
\gs (\xi) = \sum^n_{j=1} (-1)^{j+1} \xi_j\, d\xi_1 \wedge \ldots \wedge
 {d\xi_{j-1}} \wedge {d\xi_{j+1}} \wedge \ldots \wedge d\xi_n.
$$
Introducing polar coordinates, one can check that
the restriction of $\gs$ to the unit sphere $S^{n-1}$ yields 
the surface measure. 
 }%\end{extra}
\medskip

\rm We can now introduce   Wodzicki's noncommutative residue:

\def\labelenumi{(\alph{enumi})}
\def\labelenumii{(\roman{enumii})}
\begin{thm}\label{251.3}{
Let $n\ge 2$,  $A \in \Psi_{cl}(\gO)$ as in {\rm\ref{250.4}}, $x \in \gO$. 
 Denote by ${\rm tr}_E$ the trace on ${\rm Hom}(E)$ and set
$$
{\rm res}_x \, A = \left(\int_{S^{n-1}} {\rm tr_E} \,
a_{-n} (x,\xi) \gs(\xi) \right)dx_1
\wedge \ldots \wedge dx_n .
$$
This defines a density on $\gO$. Moreover, 
\begin{equation}\label{25131}
{\rm res}\, A = \int_\gO {\rm res}_x \, A 
\end{equation}
has the following properties:   
\begin{enumerate}
\item  It  only depends on the equivalence class of $A$ in $\A$. 
\item It is a trace:
${\rm res}\, [A,B] = 0$ for all $ A,B \in \A$.
\item If $\gO$ is connected, then any
other trace on $\A$ is a multiple of {\rm res}.
\end{enumerate}

{\rm
For $n=1$, the sphere $S$ consists of the two points $-1$ and $1$. 
One lets ${\rm res}\,A = \int_\gO (a_{-1}(x,-1)+a_1 (x,1))dx$. Again this is 
the unique trace on $\A$ up to multiples.

The local density
$a_{-n}(x,\xi)\sigma(\xi) \wedge dx_1\wedge\dots\wedge dx_n$  can be patched to  
a global  density $\Omega_A $ with  ${\rm res}\, A=\int_{S^*\gO}\Omega_A$: Denoting 
by   $\omega$ the  canonical symplectic form on $T^*\gO$ and 
by    $\rho $ the radial vector field one has
$$
a_{-n}\sigma\wedge dx_1\wedge\dots\wedge dx_n   =
(-1)^{n(n-1)/2}\,
\frac1{n!}(a\,\rho\,\rfloor\,\omega^n)_0,      
$$
where $(\dots)_0$ is the homogeneous component of degree 0 in an asymptotic 
expansion of $a\,\rho\,\rfloor\,\omega^n$ into homogeneous forms
        ($\rfloor$ stands for the contraction 
of  forms  with  vector fields).
}
}\end{thm}\medskip

A simple proof of Theorem \ref{251.3} was given in \cite{FGLS2}. 
Note that no continuity assumption is 
necessary for the proof of uniqueness. 
\medskip

\rem\label{local}{\rm In fact, the proof of Theorem \ref{251.3} shows even more:
We may consider the algebra of all classical symbols with support in an
open set of Euclidean space with the Leibniz product. Then any trace on
this algebra is a multiple of Wodzicki's.}

\begin{example}\label{251.4}{\rm  Let 
$A=(I-\Delta)^{-n/2}$. Then $a_{-n} (x,\xi) = |\xi|^{-n}$; hence  
$${\rm res}\, A = \int_\gO \int_{S^{n-1}}1 ~\gs(\xi) \,dx = {\rm vol}\,
S^{n-1} \cdot {\rm vol}\, \gO.$$ 
So the volume of $\gO$ can be found as a 
noncommutative residue. }\end{example}

\rm
Note that the noncommutative residue vanishes on differential operators. Also, if
the order of $A$ is $<-n$, then  $ {\rm res}\, A =0$; so res is {\em not}
an extension of  the usual operator trace. 
In fact, Connes showed that  Wodzicki's residue coincides with
Dixmier's trace on pseudodifferential operators of order $-n$ and therefore
{\em vanishes} on trace class operators. We shall see more, below.

Before, however, let us explore the relation to zeta functions and heat kernel expansions. 
\medskip

\refstepcounter{thm}

{\sc \thesection.\arabic{thm} Complex Powers.}
\label {251.6}{\rm Assume in addition that
$A$ is invertible of order $m>0$.  Then $a$ is elliptic, but we
impose a slightly stronger condition: There exists a ray 
$R_\gt = \{z\in\C: z=r e^{i\theta}, r\ge0 \}$ in $\C$ with no eigenvalue 
of the principal symbol  $a_m (x,\xi)$  on $R_\gt$  for  $  \xi \not= 0$.
The spectrum of $A$ is discrete.
Shifting $\gt$ slightly, $R_\gt$ will not intersect it. 
Seeley \cite{Seeley67} showed that:
\begin{enumerate}
\item
The norm of $(A-\gl)^{-1}$ is $O(\gl^{-1})$ on $R_\gt$, and there
exists a family of complex powers $\{ A^s: s \in \C \}$,
defined by
\begin{eqnarray*}
A^s &=& \frac{i}{2\pi} \int_{\mathcal C} \gl^s (A-\gl)^{-1} d\gl, \quad {\rm Re}\,
        s<0; \\
A^{s+k} &=& A^s A^k, \quad {\rm Re}\, s<0, k \in \N.
\end{eqnarray*}
Here
${\mathcal C} $ is the path in $\C$ going from infinity along $R_\gt$ 
to a small circle around 0, clockwise about the circle, and back along $R_\gt$.  

\item $A^s$ is a pseudodifferential operator of order $m \, {\rm Re}\, s$; $s \mapsto A^s$ is analytic.

\item
For ${\rm Re}\, s < - {n}/{m}, A^s$ is an integral operator with a
continuous integral kernel $k_s (x,y)$. For each $x \in \gO, s \mapsto k_s
(x,x)$ extends to a meromorphic map with at most simple poles  in $s_j =
\frac{j-n}{m}$, $j=0,1,\ldots$. 
There is no pole in $s=0$; the residue in $s_j$ is given by an explicit formula.
If $A$ is a {\em differential} operator, then  the residues at the positive
integers vanish.
\end{enumerate}
 
Since  $A^s$ is trace class for ${\rm Re}\, s<-n/m$ we may define the zeta function 
$$
\zeta_A (s) ={\rm trace}\, A^{-s} %= \sum \gl^{-s}_j
, ~~~~~ {\rm Re}\, s>n/m.
$$
This is a holomorphic function. It coincides with $ \int_\gO k_{-s} (x,x) dx$
hence has a meromorphic extension to $\C$ with at most simple poles in the points $s_j$.  
 
Wodzicki realized that Seeley's explicit formulas yield
\begin{equation}\label{251}
{\rm Res}_{s=-1} \zeta_A = (2\pi)^{-n}{\rm res}\, A/{{\rm ord}\, A},
\end{equation}
where ${\rm ord}\,A$ is the order of $A$, and, more generally,
\begin{equation}\label{2}
{\rm Res}_{s=s_j} \zeta_A =(2\pi)^{-n} {{\rm res}\, A^{-s_j}}/{{\rm ord}\, A}.
\end{equation}
He used this relation to define  the noncommutative residue via zeta functions:
Let $P$ be an arbitrary pseudodifferential operator. Choose $A$ as above
with ${\rm ord}\,A > {\rm ord }\, P$. 
Then also $A+uP,~u\in\R$, will meet the above requirements for small $|u|$,
 and (\ref{251}) suggests to let
\begin{equation}\label{252}
{\rm res}\, P = \frac{d}{du} {\rm res}\, (A+uP)|_{u=0} 
= (2\pi)^{n}{\rm ord}\, A\, \frac d{du}{\rm Res}_{s=-1}\, \zeta_{A+uP}(s)|_{u=0} .
\end{equation}
Moreover, Wodzicki deduced that the latter implies
\begin{equation}\label{253}
{\rm res}\, P =(2\pi)^{n}{\rm
ord}\, A\,  {\rm Res}_{s=0}\, {\rm trace}\,(PA^{-s}).
\end{equation}
}\medskip
%\end{extra}

\refstepcounter{thm}
\noindent
{\sc \thesection.\arabic{thm}  `Heat' Trace Expansions.}\label{252.5}
{In addition to  the assumptions in \ref{251.6}  we demand
that the eigenvalues
of the principal symbol matrix $a_m$ lie in a subsector of the right half-plane.
Then one can define
$$
e^{-tA} = \frac{i}{2\pi}\int_{\mathcal C'} e^{-t\gl} (A-\gl)^{-1} d\gl,
$$
where ${\mathcal C'}$ is a suitable contour around the spectrum. 
The operator $e^{-tA}$ solves the equation  $\partial_t + A =0$, which is a generalization
of the classical heat equation $\partial_t -\Delta =0$, where $\Delta$ is the Laplace-Beltrami operator.
It is easy to see that $e^{-tA}$ is trace class.
The identity 
\begin{eqnarray*}
A^{-s} = \gG(s)^{-1}\int_0^\infty t^{s-1}e^{-tA} dt
%    \int^\i_0 t^{s-1} e^{-\gl t} dt
%= \gl^{-s} \int^\i_0 (\gl t)^{s-1} e^{-\gl t} d(\gl t) 
%= \gl^{-s} \gG (s)
\end{eqnarray*}
shows that 
$\gG (s) \zeta_A (s) = \int^\i_0 t^{s-1} {\rm trace}\, (e^{-tA}) dt$
is the  Mellin transform of ${\rm trace}\,e^{-tA}$. 

The Mellin transform of a  function which is  
$\sim t^{-s_j} \ln^k t$ near $t=0$  has a 
 pole in $s_j$ of order $k+1$ and vice versa. 
Seeley's analysis of the zeta function therefore implies that 
\begin{equation}\label{eta}
{\rm trace}\, e^{-tA} \sim \sum^\i_{j=0} \ga_j (A) t^{\frac{j-n}{m}} +
\sum^\i_{k=1} (\gb_k (A) \ln t+\gb'_k(A))t^k.
\end{equation}
There is no term  $t^0\ln t$, since $\zeta_A$ is regular in $0$ while the Gamma 
function has a simple pole. For the same reason there are no terms  $t^k \ln t$ if 
$A$ is a {\em differential} operator. 

So we get 
${\rm res}\, A = (2\pi)^{n}{\rm ord}\, A \cdot \gb_1 (A)$. Moreover, we can define
the noncommutative residue for a  general pseudodifferential operator by choosing an operator
$A$ with the above properties and ${\rm ord }\,A >{\rm ord}\,P$, then letting
\begin{equation}\label{var}
{\rm res}\, P = - (2\pi)^{n}{\rm ord}\, A \frac{d}{du} \gb_1 (A+uP)|_{u=0} .
\end{equation}
Alternatively, one can establish an expansion for ${\rm trace}\,(Pe^{-tA})$ of the form 
$$ \sum^\i_{j=0} \ti \ga_j (A) t^{\frac{j-n-{\rm ord}\,P}{m}} +
\sum^\i_{k=0} (\ti \gb_k (A) \ln t+\ti\gb_k'(A))t^k$$ 
and define 
\begin{equation}\label{255}
{\rm res}\, P = - (2\pi)^{n}{\rm ord}\, A ~\ti \gb_0  .
\end{equation}
}
\medskip

\rm 
Relations \eqref{253} and \eqref{255} can be deduced rather easily using a 
symbolic calculus introduced by Grubb and Seeley \cite{GrubbSeeley93a} which relies
on pseudodifferential symbols that depend on a complex parameter in a special 
way:

\defn \label{252.7}{\rm Let $\gG$ be a sector in $\C\setminus\{0\}$. 
$S^{m,0}(\rn\times\rn;\gG)$
is the space of all functions 
$p= p(x,\xi,\mu)\in C^\i(\rn\times\rn\times\gG)$   that  are holomorphic
with respect to $\mu$ whenever $\mu\in {\rm int}\,\gG$ and  $|\xi,\mu|>\epsilon$ for some 
$\epsilon>0$, and which satisfy, for all $j\in\N_0$, 
$$\partial_z^jp(\cdot,\cdot,1/z) \in S^{m+j}(\rn\times\rn), ~~~~1/z\in\gG,$$
with uniform estimates for $|z|\le 1$ and  $1/z$ varying over a closed subcone of 
${\rm int }\, \gG$.
Closure refers to the topology of $\C\setminus \{0\}$. Finally, 
$S^{m,l}(\rn\times\rn;\gG)= \mu^lS^{m,0}(\rn\times\rn;\gG).$

A symbol $p$ in $S^{m,l}$ is said to be weakly polyhomogeneous
provided that there exists a sequence of symbols 
$p_j\in S^{m_j-l,l}$, $m_j\searrow -\infty$, $j=1,2,\ldots$, with
$p_j$ homogeneous in $(\xi,\mu)$ for $|\xi|\ge 1$ of degree $m_j$ such that
$p\sim\sum p_j$. In particular we then have $p\in  S^{m_1-l,l}$. Write
$p\in S^{m,l}_{\rm wphg}$.
}\medskip

\rm Grubb and Seeley call 
a symbol $p$  {strongly polyhomogeneous}
 of degree $m$ with respect to $(\xi,\mu)\in \rn\times (\gG\cup\{0\})$,
when $p\in C^\infty(\rn\times\rn\times(\gG\cup\{0\}))$ and there is a sequence of 
functions $p_j\in C^\infty(\rn\times\rn\times(\gG\cup\{0\}))$ that, for 
$|\xi,\mu|\geq 1,$ are
homogeneous in $(\xi, \mu)$ of degree $m-j$, with $p$ and the $p_j$ holomorphic in 
$\mu\in{\rm int}\,
\gG$ and 
$$
D^\ga_\xi D^\gb_x D^k_\mu (p-\sum_{j=0}^J p_j) = O(\n{\xi,\mu}^{m-J-|\ga|-k}),
$$
uniformly for $\mu$ in closed subsectors of $\gG\cup\{0\}$; as usual, $\langle y \rangle = (1+|y|^2)^{1/2}$ for $y \in \R^k$.  
They write $p\in S_{\rm sphg}^m $.

If $p$ is strongly polyhomogeneous of degree $m\le0$,   then
 $p\in S^{m,0}\cap S^{0,m}$ by \cite[Theorem 1.16]{GrubbSeeley93a}.

\rm An important result is the following theorem  \cite[Theorem 2.1]{GrubbSeeley93a}:

\begin{thm}\label{3.25}{
Let $p \sim \sum^\i_{j=0} p_j$ in $S^{\i,d}_{\rm wphg}$, where $p_j$ is homogeneous of degree $m_j, m_j \downarrow -\i$,
and has $\mu$-exponent $d$. Assume further that $p$ and the $p_j$ with
$m_j-d \geq -n$ are in $S^{m',d'}$ for some $m'<n'$, some $d' \in \R$.
Then $\Op p$ has a kernel $K_p (x,y,\mu)$ with an expansion along the
diagonal
$$
K_p (x,x,\mu) \sim \sum^\i_{j=0} c_j (x) \mu^{m_j+n} + \sum^\i_{k=0}
[c'_k (x) \ln \mu + c''_k (x)] \mu^{d-k}
$$
for $|\mu| \to \i$, uniformly for $\mu$ in closed subsectors of $\Gamma$.\medskip

\rm The coefficients $c_j (x)$ and $c'_{d-m_j-n} (x)$ are determined by $p_j
(x,\xi,\mu)$ for $|\xi| \geq 1$ (are ``local''), while the $c''_k(x)$ are
not in general determined by the homogeneous parts of the symbol (are
``global'').
}
\end{thm}

\rm
In addition, one easily deduces from the proof of this theorem that the
coefficient of $\mu^{d-k} \ln \mu$ is given by
\begin{equation}\label{res1}
(2\pi)^{-n}\int_{|\xi|=1} \frac{1}{k!} \partial^k_z (z^d p_j (x,\xi,\frac{1}{z}))
|_{z=0} \, \gs(\xi) ,
\end{equation}
where $j$ is such that $k=d-m_j-n$.

Only one $p_j$ contributes to this coefficient. For $d=0$,
the integrand is the coefficient of $\mu^{-k}$ in a Taylor expansion of $p_j
(x,\xi,\mu)$ into powers of $\mu^{-1}$; it is homogeneous of degree $-n$
in $\xi$.\bigskip

\refstepcounter{thm}
 {\sc \thesection.\arabic{thm} Application to Wodzicki's Residue.}\label{252.8}{ \rm
The above classes are particularly suited for the analysis of  $P(A-\gl)^{-1}$,
where $P$ is an arbitrary pseudodifferential operator and $A$ is as above, its order being
a positive integer. 
The reason is the following:
Letting $\mu^m=\gl$, the 
symbol of $(A-\gl)^{ -1} = (A-\mu^m)^{-1}$ is strongly (hence weakly) polyhomogeneous.
The composition with a $\mu$-independent symbol stays weakly polyhomogeneous
of the same order.
Applying Theorem \ref{3.25}  one obtains  an expansion
\begin{equation}
{\rm trace}\,(P(A-\gl)^{-k}) \sim\sum_{j=0}^\infty c_j
\gl^{\frac{n+{\rm ord}\,P-j}{m}-k} + \sum_{l=0}^\infty (c'_l\ln\gl + c''_l)\gl^{-l-k},\label{E1}
\end{equation}
\cite[Theorem 2.7]{GrubbSeeley93a}. Here, $k$ is chosen  larger than $(n+{\rm ord}\,P)/m$ so that the operator under
consideration is indeed trace class.
In fact, the operator $P(A-\gl)^k$ is an integral operator with a continuous kernel
$K=K(x,y,\gl)$, and the corresponding expansion holds for the values $K(x,x,\gl)$
along the diagonal.

We are interested in the coefficient of
$\gl^{-k}\ln\gl= m\mu^{-mk}\ln\mu$ in \eqref{E1}.
We just saw that we obtain the coefficient of $\mu^{-mk} \ln \mu$ by
integrating the Taylor coefficient of $\mu^{-mk}$ in the expansion of the
symbol of $P(A-{\mu^m})^{-k}$, more precisely the component  of homogeneity
$-n$. This is the highest power of $\mu$ that can possibly occur. The
asymptotic expansion formulae show that it only is involved in the term
$$[p(x,\xi) (a(x,\xi)-\mu^m)^{-k}]_{-n}\,.$$
Writing
$$(a-{\mu^m})^{-k}=(-1)^k\mu^{-mk}(1-{a}/{\mu^m})^{-k}=(-1)^{k}\mu^{-mk}\left(1+\ldots(
{\rm powers\ of\ } a/\mu^m\right))$$
we see that the result is independent of $a$ and that the coefficient of
$\gl^{-k}\ln \gl$ is
\begin{equation}
\label{res2}
(2\pi)^{-n}\frac{(-1)^k}{{\rm ord}(A)}\, \int_\Omega \int_{|\xi|=1}
p_{-n}(x,\xi)\gs(\xi)\,dx.
\end{equation}

One next notes the identities
\begin{eqnarray}
\label{AS}
A^{-s} &= &
\frac{ik!}{2\pi(s-1)\ldots(s-k)} \int_{\mathcal C} \gl^{k-s}(A-\gl)^{-k-1}d\gl \\
&=& \gG(s)^{-1}\int_0^\infty t^{s-1}e^{-tA} dt\nonumber\\
{}e^{-tA}&=&\frac{ik!} {2\pi}(-t)^{-k}\int _{\mathcal C'} e^{-t\gl} (A-\gl)^{-k-1}
d\gl.\label{ETA}
%= \frac1{2\pi i}
\end{eqnarray}
They are obtained from the usual Dunford integrals via integration by parts and the
relation
$$\partial^k_{\gl}(A-\gl)^{-1} = k!~(A-\gl)^{-k-1}.$$

They imply that
\begin{eqnarray}
{\rm trace}\,(Pe^{-tA}) &\sim&\sum_{j=0}^\infty \ti c_j t^\frac{j-n-{\rm ord}\,P}{m} + \sum_{l=0}^\infty (\ti c'_l\ln t + \ti c''_l)t^l,\label{E2}\\
\Gamma(s) {\rm trace}\,(PA^{-s}) &\sim&\sum_{j=0}^\infty \frac{\ti c_j}{s+\frac{j-n-{\rm ord}\,P}{m}} + \sum_{l=0}^\infty \left(\frac{-\ti c'_l}{(s+l)^{2}} +\frac{ \ti c''_l}{s+l}\right)\label{E3}.
\end{eqnarray}
It is  no longer necessary that the order of $A$ be large.
In \eqref{E3}, the tilde indicates that the left hand side is meromorphic with
poles as indicated by the right hand side. The coefficients $\ti c_j$, $\ti c'_l$, and $\ti c''_l$
are multiples of the corresponding $ c_j$, $c'_l$, and $c''_l$, the factors are universal constants independent of $A$ and $P$.
In particular, we see from (\ref{res2}), (\ref{AS}), and (\ref{E3}) that the coefficient
$\ti c_0'$ is given by
\begin{equation}\label{c}
\ti c_0' = -
\frac{ (2\pi)^{-n}}{{\rm ord}\, (A)}
\int_\Omega \int_S p_{-n} (x,\xi) \, \gs(\xi) \, dx
= -\frac{ (2\pi)^{-n}}{{\rm ord}\, (A)}\, {\rm res}\, (P).
\end{equation}
 }\bigskip

\centerline{\sc Dixmier's Trace}\medskip

\rm
Let us now turn to the analysis of Dixmier's trace.
By ${\mathcal K}(H)$ denote the ideal of compact operators on the
Hilbert space $H$. In view of Example \ref{250.2}
it had been an open question, whether every completely additive trace is
proportional to the standard trace on the subset of ${\mathcal K}(H)$ where it is finite.
Dixmier \cite{Dixmier66} showed that the answer is `no'  by explicitly constructing counter-examples.
We give a short review,   following
Connes \cite{Connes94} in presentation and terminology.\bigskip

\refstepcounter{thm}
 {\sc \thesection.\arabic{thm}
The Spaces  ${\mathcal L}^{(1,\i)} (H)$ and ${\mathcal L}_0^{(1,\i)} (H)$. }\label{9.1}{\rm
Let $H$  be an infinite-dimensional Hilbert space,
$T \in {\mathcal K}(H)$, and $|T| = (T^* T)^{1/2}$. Let
$ \mu_0(T)\ge \mu_1(T) \ge \ldots  $
be the sequence of the singular values of $T$, i.e. eigenvalues of $|T|$,
repeated according to their multiplicity.
%It is well-known that
%\begin{equation}
%\hspace*{-1em}\mu_j(T) = \inf \{ \|T-F\|: {\rm rank }\, F = j\}=
%\min\{\|T|_{E^\bot}\|: \dim E = j\}.
%\end{equation}
%We deduce from (1.9) and (1.10) that
%\begin{eqnarray}
%{\rm Res}_{s=0} \mbox{trace}\, (PA^{-s} = (2\pi)^{-n}{\rm ord}A
%\int_\Omega\,\int_{|\xi|=1}p_{-n}(x,\xi) d\sigma (\xi)dx
%\end{eqnarray}

We define $\gs_N (T) = \sum^N_{j=0} \mu_j (T)$ and let
$${\mathcal L}^{(1,\i)} (H) = \{ T \in {\mathcal K} (H): \gs_N (T) = O(\ln N) \},$$
endowed with the
norm
$$
\| T \|_{1,\i} = \sup_{N \geq 2} \frac{\gs_N (T)}{\ln N}.
$$
%We have a natural subspace
%${\mathcal L}_0^{(1,\i)} (H) = \{ T \in {\mathcal K} (H): \gs_N (T) = o (\ln N) \}$.
%\subseteq \L^{(1,\i)} (H)$.
%${\mathcal L}^{(1,\i)}_0,
%$ {\mathcal L}^{(1,\i)}$
This clearly is a two-sided ideal in ${\mathcal L}(H)$.
}\medskip

\refstepcounter{thm}
 {\sc \thesection.\arabic{thm} Ces\`aro Mean. }\label{9.4}{\rm
We define the Ces\`aro mean $Mf $ for $f\in L^\i (1,\i) $ by
%$\to C_b (1,\i)$ by
$$
(Mf)(t) = \frac{1}{\ln t} \int^t_1 f(s) \frac{ds}{s}.
$$
The function $Mf$ is continuous and bounded;  $M: L^\i (1,\i) \to C_b (1,\i)$
is continuous.  Moreover,  $M1=1$ and, for $\gl>0$, 
$M(f(\gl\cdot)) - Mf \in C_{b(0)} (1,\i)$.
The subscript $(0)$ indicates that the function vanishes at infinity.
}\medskip

\refstepcounter{thm}
 {\sc \thesection.\arabic{thm} The `limit' $\lim_\go$ and Dixmier's trace. }\label{9.5}{\rm% %Let 
%$T_1,T_2 \in {\mathcal L}^{(1,\i)}$ be positive and
%$$
%\ga_N = \frac{\gs_N (T_1)}{\ln N}, ~~ \gb_N = \frac{\gs_N (T_2)}{\ln N}, ~~
%\gg_N = \frac{\gs_N (T_1+T_2)}{\ln N} .
%$$
%Then $ \{ \ga_N\}, \{ \gb_N \}$, and $ \{ \gg_N \}$ are bounded sequences. By \ref{9.3}
%we have 
%\begin{equation}
%\gg_N \le \ga_N+ \gb_N \leq ({\ln 2N}/{\ln N})\, \gg_{2N},
%\end{equation}
%but in general no convergence. 
We  embed $\cL^\i (H)$ into $L^\i (1,\i)$
  by associating to the sequence
$\{ a_j \}$ the function $ f_{ \{ a_j \} }$
which has the value $a_j$  on the interval $[j,j+1[$, $j=1, 2,\ldots$.
Next we choose a linear form $\go$ on $C_b (1,\i)$ with
(i) $\go \geq 0$, (ii) $ \go(1) =1$, and  (iii) $\go(f) =0$ for $ f \in C_{b(0)}$. Then we 
define $\lim_\go \{ a_j \} = \go (M f_{\{a_j\}})$ with the help of Ces\`aro's mean. 

Note that  $\lim_\go$ coincides with the usual
limit on convergent sequences by  (ii) and  (iii).
%Furthermore, $\lim_\go a_{2N} = \lim_\go a_N$.}

{For  a positive operator $T\in \cL^{(1,\i)} (H)$ we now let
$$
{\rm Tr}_\go (T) = \lim_\go \frac{1}{\ln N} \sum^N_{n=0} \mu_n (T) .
$$
As Proposition \ref{9.7}(a), below,  states, ${\rm Tr}_\go$ is additive on positive $T$'s. We can therefore
extend it uniquely to a linear map on $\cL^{(1,\i)} (H)$, also denoted ${\rm Tr}_\go$. 
}\medskip

\begin{prop}\label{9.7}{Let 
$T,T_1,T_2 \in \cL^{(1,\i)} (H), S \in \cL(H)$.
\bli
\item ${\rm Tr}_\go (T_1+T_2) = {\rm Tr}_\go (T_1) +
 {\rm Tr}_\go (T_2)$ for positive $T_1,T_2 $.
\item ${\rm Tr}_\go (T) \geq 0$ if $T\geq0$.
\item If $S$ is  invertible, then  ${\rm Tr}_\go (STS^{-1}) = {\rm Tr}_\go (T)$.
In particular, ${\rm Tr}_\go$ is independent of the inner product in $H$.
\item ${\rm Tr}_\go (ST) = {\rm Tr}_\go (TS)$.
\item ${\rm Tr}_\go$ 
% \equiv 0$ on $\cL^{(1,\i)}_0$, so it 
vanishes on  trace class operators.
\eli
}\end{prop}

%\noindent\Proof.  (a) follows from \ref{9.5}(1) together with
%the last remark in \ref{9.5}.
%We only have to check  (c) for positive $T$. Then use  \ref{9.3}(a)/(b).
%Finally (c) implies  (d) first for invertible $S$, then for arbitrary
%$S$ by adding a large multiple of the identity and using linearity.
%\eproof

\begin{example}\label{9.8}{\rm Consider the operator
$(1-\Delta)^{-n/2} : L^2 (\T^n) \to L^2 (\T^n)$,
where $\Delta $ is the Laplacian.
The eigenvalues of   $\Delta$ are known to be the lengths
$ | k |^2$ as $ k$ varies over $\Z^n$, so the
eigenvalues of  $(1-\Delta)^{-n/2}$ are $ (1-|k|^2)^{-n/2}$. 

Let us show that  $(1-\Delta)^{-n/2} \in \cL^{(1,\i)}$
and ${\rm Tr}_\go (1-\Delta)^{-n/2} =
{\gO_n}/{n}$, independent of $\go$ with
$\gO_n = {\rm vol}\, S^{n-1}$:
We let $N_R$ denote the number of
lattice points in $B_R$, the ball of radius $R$. Clearly, 
$N_R\sim {\rm vol}\,B_R$, hence $\ln N_R \sim n\ln R$. Moreover,
\begin{eqnarray*}
\sum_{|k|\le R} (1+|k|)^{n/2} &\sim& \gO_n \int_0^R (1+r^2)^{-n/2}r^{n-1}dr\\
&\sim& \gO_n \int_1^R r^{-1}dr=\gO_n\ln R.
\end{eqnarray*}
We conclude that 
$$
(\ln N_R)^{-1} \sum_{|k|\le R} (1+|k|)^{n/2} \sim \frac{\gO_n\ln R}{n\ln R}= \frac{\gO_n}{n}.
$$

Recall that
for $(1-\Delta)^{-n/2}: L^2 (\T^n) \to L^2 (\T^n)$ we had computed in
\ref{9.7} that
$$
{\rm res}\, (1-\Delta)^{-n/2} = {\rm vol}\, S^{n-1} {\rm vol}\, \T^n =
\gO_n (2\pi)^n.
$$
}\end{example}

%\subsection*{Connes' Theorem}

\noindent
In one special case we therefore have checked the following result:
\begin{thm}\label{10.1}{ \sc (Connes 1988) }{Let 
%$M$ be closed, compact, $n$-dimensional, $E$ a vector bundle over $\gO$, and
%$P: L^2 (M,E) \to L^2 (M,E)$
$P \in \Psi_{cl} (\Omega)$ be a pseudodifferential operator of order $-n$. Then
\bli
\item $P \in \cL^{(1,\i)} (L^2(\Omega,E))$.
\item $    {\rm res}\, P= (2\pi)^n n{\rm Tr}_\go P$, and the right hand side is 
independent of $\go$.
\eli
}\end{thm}

Proofs can be found in Connes' original paper \cite{Connes88} as well as
in  \cite{Varilly93} and \cite{Schrohe96c}.

\section{Boutet de Monvel's Calculus}
In the following, $X$ is an $n$-dimensional manifold with boundary
$\partial X$ embedded in $\gO$. We assume that $\dim X>1$.

\begin{defn}\label{1.4}{\rm 
\bli 
\item   $\cS (\R^n_\pm)$ denotes the rapidly decreasing functions on
      $\R^n_\pm$, i.e., $\cS (\R^n)|_{\R^n_\pm}$.
\item  $e^\pm$ is the operator of extension (by zero) of
      functions on $\R^n_\pm$ to functions on $\R^n$, while $r^+$ is the
restriction operator  from $\R^n$ to $\R^n_+$. 
We also write $e^+$  for extension by zero from $X$ to $\gO$ and $r^+$
for restriction from $\gO$ to $X$.

\item $H^s (X)$ is the Sobolev space of order $s$ on $X$, i.e. the
      restrictions of distributions in $H^s (\gO)$ to $X$, while
     $H^s (\partial X)$ is the Sobolev space at the boundary.

\item $
      H^+ = \{ (e^+ u)^\wedge: u \in \cS (\R_+) \}; \quad H^-_0 = \{ (e^-
      u)^\wedge: u \in \cS (\R_-) \} .
      $
      The hat  denotes the Fourier transform: 
$$\hat f(\tau) = (2\pi)^{-1/2}\int e^{-it\tau} f(t)dt.$$
$H'$ is the space of all 
      polynomials. We let $H=H^+ \oplus H^-_0 \oplus H'$.  The functions in $H^+ \oplus H^-_0$ satisfy $h(z) = O(\langle z \rangle^{-1})$ as $z \to \infty$. We let $H^-_d = \{ h \in H^-_0 \oplus H': h(z) = O(\langle z \rangle^{d-1}) \}$. 
 
\item We define the operator $\Pi':H\to \C$ as follows: 
Given a function $h=h_+\oplus h_-\oplus p$ in $H$ with $h_+ = (e^+u)^\wedge$
we let $\Pi'h = (2\pi)^{-1/2} u(0).$ It is easy to see that, for $h\in H\cap L^1$, 
\begin{equation}\label{Pi}
\Pi'h = \frac 1{2\pi} \int h(\tau)d\tau.
\end{equation}
\eli
}
\end{defn}

\refstepcounter{thm}
 {\sc \thesection.\arabic{thm} Boutet de Monvel's Algebra. }\label{1.6}{\rm 
For detailed introductions see Boutet de Monvel \cite{BdM},
Rempel-Schulze \cite{RS}, or Grubb \cite{Grubb86}. A short version can be
found in Schrohe and Schulze \cite{SchroheSchulzeI}.

We consider matrices of operators acting on sections of vector bundles
$E_1,E_2$ over $X$ and $F_1,F_2$ over $\partial X$. An operator of order
$m \in \Z$ and type $d$ is a matrix
$$
A = \left( \begin{array}{cc} P_+ + G & K \\ T & S \end{array} \right) :
\begin{array}{ccc} C^\i (X,E_1) && C^\i (X,E_2) \\ \oplus & \to &\oplus \\
C^\i (\partial X, F_1) && C^\i (\partial X, F_2) \end{array} .
$$
%In view of the fact that we shall deal with endomorphisms we assume
%$E_1=E_2=E$, ~$F_1=F_2=F$. 
 ${\mathcal B}^{m,d} (X)$ denotes  the collection of all operators
of order $m$ and type $d$, and $\B$ is the union over all $m$ and $d$.

{\em Convention}:
In order to avoid superfluous notation, we
shall omit ``classical''. All operators, however, are assumed
to be classical.

Let me now give a short account of the various entries.
\begin{itemize}
\item $P$ is a classical pseudodifferential operator of
order $m$ on $\gO$. The subscript ``$+$'' indicates that we consider $P_+ =
r^+ P e^+$. The operator $P$ is supposed to have the
transmission property; this means that, for all $j,k,\ga$,
the homogeneous component $p_j$
of order $j$ in the asymptotic expansion  of the symbol $p$ of $P$
% \sim \sum^\i_{k=0} p_{m-k}$
%of the complete symbol $p=p(x',x_n,\xi',\xi_n)$ of $P$
in local coordinates near the boundary satisfies
$$
\partial^k_{x_n} \partial^\ga_{\xi'} p_j (x',0,0,+1) = (-1)^{j-|\ga|}
\partial^k_{x_n} \partial^\ga_{\xi'} p_j (x',0,0,-1) .
$$
\item $G$ is a  singular Green operator (s.G.o.) of order $m$
and type $d$. Recall that a singular Green symbol of order $m$ and
type $d$ is a function
$$
g = g (x',\xi', \xi_n, \eta_n) {\rm \ on\ } \R^{n-1} \times \R^{n-1} \times \R
\times \R
$$
such that
\begin{equation}\label{g}
g (x', \xi', \n{\xi'} \xi_n, \n{\xi'} \eta_n) \in
S^{m-1}(\R^{n-1} \times \R^{n-1}) 
\hat{\otimes}_\pi H^+ \hat{\otimes}_\pi H^-_d.
\end{equation}
Being classical, it has an asymptotic expansion
$$
g \sim \sum_{k=1}^\infty g_{m-k},
$$
where $g_j$ is homogeneous of degree $j$ in $(\xi', \xi_n, \eta_n)$ for $ |\xi'| 
\geq 1$.

Such a singular Green symbol introduces an operator ${\rm op}_Gg$ on
${\mathcal S}(\R^n_+)$ by
\begin{eqnarray}   \label{opg}
\lefteqn{(\Op_G g) u (x', x_{n})}\\
&=& (2 \pi)^{-\frac{n-1}{2}}
  \int e^{ix \xi}~ \Pi'_{\eta_n}\left(g (x', \xi',\xi_n, \eta_n)
   (e^+ u) ^\wedge (\xi', \eta_n) \right) d \xi.\nonumber
\end{eqnarray}
The subscript $\eta_n$ only indicates that $\Pi'$ acts with respect to this
variable.
This reduces to $(2 \pi)^{-\frac{n+1}{2}} \iint e^{ix \xi} ~g (x', \xi',\xi_n, \eta_n)
   (e^+ u) ^\wedge (\xi', \eta_n) d \eta_n d \xi$ if $d=0$.
%There is a certain freedom in the choice of the power of $2\pi$; we here follow
%\cite[Equation 2.3.2.1(5)]{RS}.
A s.G.o.~of order $m$ and type $d$ is an operator on $C^\infty(X)$
which is locally given by such an $\Op_G g$.
\item $T$ is a trace operator of order $m$ and type $d$.
Locally, it is given by a type $d$
trace symbol, i.e., a smooth function $t = t(x', \xi', \xi_n)$ such
that $t(x', \xi', \n{\xi'} \xi_n) \in S^{m -1/2} (\R^{n-1} \times
\R^{n-1} )\hat{\otimes}_\pi H^-_d$. The symbol $t$  induces a trace
operator
${\Op}_{T} t:{\mathcal S}(\R^n_+) \rightarrow {\mathcal S}(\R^{n-1})$ by
$$
({\Op}_{T} t) u (x) =
(2\pi)^{- \frac {n-1}{2}} \int e^{ix'\xi'} \Pi'_{\xi_n} (t (x',\xi', \xi_n)
(e^+u)^\wedge (\xi', \xi_n) ) d \xi'.$$
The classicality of $t$ here refers to the fact that $t$ has an expansion
$$
t \sim \sum_{k=0}^\infty t_{m-k-\frac{1}{2}}
$$
where $t_j$ is homogenous of degree $j$ in $(\xi', \xi_n)$ for $|\xi'| \geq 1$.

\item $K$ is a potential operator (or Poisson operator) 
of order $m$. It is given locally in the form
${\Op}_{K} k$, where $k = k(x', \xi', \xi_n)$ is a potential symbol of
order $m$,
i.e. satisfies $k(x', \xi', \n{\xi'} \xi_n)
\in S^{m -1/2} (\R^{n-1} \times \R^{n-1}) \hat{\otimes}_\pi H^+$,  and where
${\Op}_{K} k : {\mathcal S} (\R^{n-1}) \rightarrow {\mathcal
S} (\R_+^n)$ is given by
$$
({\Op}_{K} k) u (x', x_n) = (2\pi)^{-\frac{n-1}{2}} \iint e^{ix\xi} k 
(x', \xi', \xi_n) \hat{u} (\xi') d \xi_n d \xi'.
$$
Similarly as above we have an asymptotic expansion $k\sim\sum_{j=0}^\infty
k_{m-j-1/2}$ with $k_j$ homogeneous of degree $j$ in $(\xi',\xi_n)$
for $|\xi'|\ge 1$.

\item Finally, $S$ is a pseudodifferential operator of order $m$ along the 
boundary.
\end{itemize} \medskip
{\em Regularizing Elements.} The intersection over the 
space  of s.G.o. of all orders $m$ is called the space of regularizing s.G.o.s (the type may be arbitrary). 
The regularizing s.G.o.s of type zero
are integral operators with a smooth integral kernel.
Similarly we obtain the notions of regularizing trace and potential operators.
The space of matrices $A$ as above with all entries regularizing is denoted 
by $\B^{-\infty}$. 

}\medskip

\refstepcounter{thm}
 {\sc \thesection.\arabic{thm} Basic Properties. }\label{1.7}\rm 
We shall rely on  a few facts.
\bli
\item The sum in the upper left corner is direct up to regularizing
pseudodifferential operators. In other words,
if $P_+=G$ for a pseudodifferential operator 
$P$ and a singular Green operator $G$, then $P$ is regularizing. Indeed this
is a consequence of the following:
\item Given a function $\omega$ which is equal to one near the boundary
and a s.G.o.~$G$ 
of order $m$ and type $d$, the composition $(1-\omega)G$ is a regularizing
s.G.o.~of type $d$, while $G(1-\omega)$ is regularizing of type zero.
\item ${\mathcal B}^{m,d}$ is a Fr\'echet space. The composition of the above
operator matrices yields a continuous map
\begin{equation}\label{compos}
{\mathcal B}^{m,d} \times {\mathcal B}^{m',d'} \rightarrow {\mathcal B}^{m+m',
\max\{m'+d,d'\}}
\end{equation}

Here we assume that the bundles the operators act on are such that the
composition makes sense. In order to illustrate this somewhat more write
$$
A = \left( \begin{array}{cc} P_+ + G & K \\ T & S \end{array} \right) \in 
{\mathcal B}^{m,d}, \quad
A' \in  \left( \begin{array}{cc} P'_+ + G' & K' \\ T'& S' \end{array} \right) \in 
{\mathcal B}^{m',d'}
$$

The composition $AA'$ is obtained by multiplication of the matrices. The 
resulting compositions of the entries then are of the expected kind; the order 
of all entries is $m+m'$.
$TK'$, for example, is a pseudodifferential operator
on the boundary; $KT'$, $P_+G'$, and $GG'$ are
s.G.o.s~of type $d'$, while $GP'_+$ is a s.G.o.~of type $m'+d$. For the
composition of the pseudodifferential parts we have
$$
P_+ P'_+ = (PP')_+ + L(P, P').
$$
Here $PP'$ is the usual composition of pseudodifferential operators, while 
$L(P,P')$, the so-called {\em leftover term}, is a s.G.o.~of type $m'+d$.

\item It follows that the space $\B_0^{m,d}$
of all operator matrices,
where the pseudodifferential part $P_+$ is zero, forms a two-sided ideal
in the sense that \eqref{compos} restricts to mappings
\begin{eqnarray*}
{\mathcal B}_0^{m,d} \times {\mathcal B}^{m',d'} &\rightarrow& {\mathcal B}_0^{m+m',
\max\{m'+d,d'\}}\mbox{\rm ~~and}\\
{\mathcal B}^{m,d} \times {\mathcal B}_0^{m',d'} &\rightarrow& {\mathcal B}_0^{m+m',
\max\{m'+d,d'\}}.
\end{eqnarray*}

\item The linear span of all compositions $KT'$ where $K$ is a potential
operator of order $m$ and $T'$ is a trace operator of order $m'$ and type
$d'$ is dense in the space of all s.G.o.s of order $m+m'$ and type $d'$.
The compositions $TK'$ similarly span the space of pseudodifferential symbols 
of order $m+m'$. 

\item By (\ref{opg}), a singular Green symbol $g$ of type zero defines, for fixed $x',\xi'$, an operator
$g(x',\xi',D_n): {\mathcal S}(\rp) \to {\mathcal S}(\rp)$ by 
$$g(x',\xi',D_n)u(x_n) = 
\frac 1{2\pi}
\int e^{ix_n\xi_n}g(x',\xi',\xi_n,\eta_n)\hat u(\eta_n) d\eta_n.$$
It  extends to a trace class operator on $L^2(\rp)$; the trace is given by
\begin{eqnarray*}
%{\rm tr}\,g(x',\xi')& =&
\lefteqn{{\rm trace}\, g(x',\xi', D_n) }\\
&=& (2\pi)^{-{1/2}}\int g(x',\xi',\xi_n,\xi_n)d\xi_n = (2\pi)^{1/2}\,\Pi'_{\xi_n}g(x',\xi',\xi_n,\xi_n),
\end{eqnarray*}
where the subscript $\xi_n$ indicates the variable $\Pi'$ is acting on.
For convenience we slightly change the norming  factor and define the functional
\begin{equation}\label{tr}
{\rm tr}\,g(x',\xi')= \Pi'_{\xi_n}g(x',\xi',\xi_n,\xi_n).
\end{equation}
In this way it even extends to the case where $g$ is of type $d$.
%$g=\sum_{j=0}^d g_j\xi_n$ with 
%$g_j$ of type zero (corresponding to 
%$G$ of type $d$ by \eqref{G}). 

Note: If $g(x',\xi',\xi_n,\eta_n)$ is homogeneous of degree $j$ in $(\xi',\xi_n,\eta_n)$ 
for $|\xi'|\ge 1$, then ${\rm tr}\,g(x',\xi')$ is homogenous in $\xi'$ of degree $j+1$.
\item An operator $A\in {\mathcal B}^{m,d}$ extends to a bounded operator
$$A:\vectd{H^s(X,E_1)}{\oplus}{H^s(\partial X,F_1)}\to
\vectd{H^{s-m}(X,E_2)}{\oplus}{H^{s-m}(\partial X,F_2)}$$
for each $s>d-1/2$.

\eli
\medskip

\section{The Noncommutative Residue for Manifolds with Boundary}

We use the notation of the previous section. In view of the fact that
we are dealing with endomorphisms, we assume that $E_1 = E_2 = E$ and
$F_1=F_2=F$.

It was already mentioned in the introduction that it is not
sufficient to consider merely the algebra of usual  classical
pseudodifferential symbols (with the Leibniz product as composition),
if one intends to find a trace functional
for boundary value problems. 
Wodzicki made the following observation:

\lem\label{notrace}{There is no non-zero trace
on $\Psi_{cl}(\gO)/\Psi^{-\infty}(\gO)$ whenever $\gO$ is noncompact or
has a boundary.}\medskip

\rm {\Proof.}  We know from Remark \ref{local} that the only candidate for
a trace is a multiple of the noncommutative residue.
For noncompact manifolds, however, formula \eqref{25131} does not make
sense.

Next suppose that $\gO$ has a boundary and there exists a trace. Consider the symbols with support in
$U\times\rn$, where $U$ is a neighborhood of the boundary
of the form $U=U_0\times[0,1)$ with $U_0$ open in $\R^{n-1}$.
Every trace necessarily vanishes on symbols $p=p(x,\xi)$ of the form
$p=\partial_{x_k}q$ for another symbol $q$, since
$\partial_{x_k}q$ is the symbol of the commutator
$[\partial_{x_k},\Op q]$. On the other hand,
each symbol with support in $U\times\rn$
has a primitive with respect to $x_n$.
Thus it is a derivative and hence has trace zero.
The trace therefore vanishes.\eproof\medskip

For more details see \cite[Theorem 3.2]{FGLS2}.\medskip

Note, in particular, that the expression
$$
\tau (P) = \int_X \int_S {\rm tr}_E p_{-n} (x,\xi) \gs(\xi) dx
$$
does {\em not} define a trace if we admit pseudodifferential operators
of all orders (it does define a trace on the algebra of operators of order $-n$).\bigskip

On a manifold with boundary, however, it is much more natural to work with
Boutet de Monvel's algebra. There, one obtains the following
theorem \cite{FGLS1}, \cite{FGLS2}.

%\begin{dedef\label{}{}}

\begin{thm}\label{FGLS}{Let $\Omega$ be connected of dimension $\ge 2$, let
$$
A = \left( \begin{array}{cc} P_+ + G & K \\ T & S \end{array} \right) \in
{\mathcal B},
$$and denote by $p,g,$ and $s$ the local symbols of $P,G,$ and $s$
respectively. Define
\begin{eqnarray}
{\rm res}\, A &=& \int_X \int_S {\rm tr}_E p_{-n} (x,\xi) \sigma(\xi) dx\nonumber\\
&& {} + 2 \pi \int_{\partial X} \int_{S'} \left\{ {\rm tr}_E ({\rm tr}\,
g _{-n}) (x',\xi') + {\rm tr}_F s_{1-n}(x',\xi') \right\}
\sigma'(\xi') dx' .\label{fgls}
\end{eqnarray}
Then
\bli
\item ${\rm res}$ is a trace: ${\rm res} [A,B] = 0$ for all $A,B \in {\mathcal
B}$.
\item It is the unique continuous trace on ${\mathcal B}/{\mathcal B}^{-\infty}$.
\eli

\rm Here, $\sigma'$ is the analog of the form $\sigma$ in the $(n-1)$-dimensional case (for
$n=2$, we replace the integral $\int_{S'} f(x',\xi') \sigma'(\xi')$ of a function $f$ over
$S^{n-1}$  by  $f(x',1)+f(x',-1)$). The subscripts $-n$ and $1-n$ indicate that
we are considering the components of the corresponding degree of homogeneity.
Note that ${\rm tr}$ is the trace on singular Green symbols introduced in \eqref{tr}
and that ${\rm tr}\, g_{-n}$ is homogeneous of degree $1-n$ in $\xi'$.

 This trace `res' reduces to Wodzicki's noncommutative residue
if $\partial X= \emptyset$ and was therefore called the
noncommutative residue for manifolds with boundary
in \cite{FGLS1}, \cite{FGLS2}.\footnote{The sign between the two terms in
\eqref{fgls} should indeed be `$+$', not `$-$'.}}
\end{thm}

A proof can be found in \cite{FGLS2}. Let me give a sketch of at least part of the 
argument: Using \ref{1.7}(a), the algebra $\B/\B^{-\infty}$ can be decomposed 
$$\B/\B^{-\infty}= \Psi_{\rm tr}(X)/\Psi^{-\infty}(X) \oplus \B_0/\B^{-\infty}
$$
into the algebra 
$\Psi_{\rm tr}(X)/\Psi^{-\infty}(X)$  of all classical pseudodifferential
symbols with the transmission property, and the quotient of $\B_0$ -- i.e.~the matrices in Boutet de Monvel's calculus that have zero pseudodifferential part -- modulo the regularizing elements. 

On  symbols supported in the interior, any trace will be a
multiple of the noncommutative residue. This suggests that the
trace is of the above form on the pseudodifferential part. 

Next, on the algebra $\B_0/\B^{-\infty}$, one may apply the following lemma: 

\begin{lem}\label{eric}
Let $\tau:\,\B_0/\B^{-\infty}\to\C$ be linear. Then $\tau$ is a trace on $\B_0/\B^{-\infty}$
if and only if there exist a trace $\tau_1$ on
the algebra of singular Green symbols and a trace $\tau_2$ 
on $\Psi(\partial M)/\Psi^{-\infty}(\partial M)$ such that
for all trace symbols $t$ and all potential symbols $k$
\begin{eqnarray}
\tau_1(k\circ t)=\tau_2(t\circ k)
\end{eqnarray}
and
\begin{eqnarray}
\tau\left\{\left(\begin{array}{cc}
g & k\\
t & s
\end{array}\right)\right\}=\tau_1(g)+\tau_2(s)\,.
\end{eqnarray} 
\end{lem}

\Proof. 
Given   $\tau$, define $\tau_1(g) = \tau
\left\{\left(\begin{array}{cc}
g & 0\\
0 & 0
\end{array}\right)\right\}
$
and
$\tau_2(s) =  \tau
\left\{\left(\begin{array}{cc}
0 & 0\\
0 & s
\end{array}\right)\right\}
$.
The identities

\begin{eqnarray}
0 &=&\label{c1}
\tau \left\{
\left[\left(\begin{array}{cc}
0 & 0\\
0 & 1
\end{array}\right),
\left(\begin{array}{cc}
0 & 0\\
t & 0
\end{array}\right)\right]\right\}
=\tau \left\{
\left(\begin{array}{cc}
0 & 0\\
t & 0
\end{array}\right)
\right\};\\
0 &=&\label{c2}
\tau \left\{
\left[\left(\begin{array}{cc}
0 & k\\
0 & 0
\end{array}\right),
\left(\begin{array}{cc}
0 & 0\\
0 & 1
\end{array}\right)\right]\right\}=
\tau \left\{
\left(\begin{array}{cc}
0 & k\\
0 & 0
\end{array}\right)
\right\};\\
0 &=&\label{c3}
\tau \left\{
\left[\left(\begin{array}{cc}
0 & k\\
0 & 0
\end{array}\right),
\left(\begin{array}{cc}
0 & 0\\
t & 0
\end{array}\right)\right]\right\}=
\tau\left\{
\left(\begin{array}{cc}
k\circ t & 0\\
0& -t\circ k
\end{array}\right) \right\}
\end{eqnarray}
imply the assertion. Conversely, given traces $\tau_1$ and $\tau_2$ satisfying the 
compatibility condition, it is easily checked that the above $\tau$ 
defines a trace.  \eproof\bigskip

\cor\label{sgtrace}{\rm  We know that the zero map is a trace and  that the span of the $k\circ t$ is dense in the space of all singular Green symbols while the span of the 
$t\circ k$ are the pseudodifferential symbols. Hence every continuous trace $\tau_1$
on the singular Green symbols induces a trace on $\Psi_{cl}(\partial 
X)/\Psi^{-\infty}(\partial X)$.
The latter must be a  multiple of Wodzicki's residues on the boundary, 
${\rm res}_{\partial X}$. 

So we do not have much choice: 
Suppose we have a trace $\tau$ and that $\tau_2 = c\,{\rm res}_{\partial X}$. 
 According to the above formulae,  the composition  $k \circ t$ has the singular
Green symbol $g(\xi_n,\eta_n) =  k(\xi_n)t(\eta_n)$. We conclude that
$$\tau_1 (g) =\tau_1(k\circ t)
= c\, {\rm res} _{\partial X}(t\circ k) = c\,{\rm  res}_{\partial
X} \Pi'(kt) =c\, {\rm res}_{\partial X}{\rm tr}\,g(x',\xi')$$
with the trace  ${\rm tr}$ introduced in \eqref{tr}.
For  simplicity we have ignored the  $(x',\xi')$ variables in the composition --
they  only cause notational complications.

On the subalgebra $\B_0/\B^{-\infty}$ we therefore have precisely one continuous  trace up to multiples, namely
\begin{eqnarray*}
{\rm res}\, \left(\begin{array}{cc}
G & K\\
T & S
\end{array}\right) =  
  \int_{\partial X} \int_{S'} \left\{ {\rm tr}_E ({\rm tr}\, 
g _{-n}) (x',\xi') + {\rm tr}_F s_{1-n}(x',\xi') \right\}
\sigma'(\xi') dx' .
\end{eqnarray*}
}\bigskip

\rm These two simple considerations show that any candidate for a trace must be of the form 
in Theorem \ref{FGLS}. Yet it is far from obvious that \eqref{fgls} indeed defines 
a trace; in fact, this requires rather delicate computations, cf.~\cite{FGLS2}. 

\section{Dixmier's Trace for Boundary Problems}
\rm We shall now analyze the connection between the noncommutative residue and Dixmier's trace. The exposition follows the original paper, Nest and Schrohe  \cite{NestSchrohe98}.
We first show that there is indeed a subclass of Boutet de Monvel's algebra 
on which Dixmier's trace makes sense. 

Like in the previous section, the operators are supposed to be endomorphisms, hence $E_1=E_2=E$ and $F_1=F_2=F$. 
We shall work on the Hilbert space  $H = L^2 (X,E) \oplus L^2 (\partial X,F)$. 
Operators of type $d>0$ can not be defined on this space in general, we shall therefore
work with type zero operators. 

It is well-known that a bounded  operator $L^2(\partial X) \to H^{n-1}(\partial X)$
defines an element of $\cL^{(1,\i)}(L^2 (\partial X))$. This operator is even trace class on
$L^2(\partial X)$,  if its range   is contained in  
$H^{n}(\partial X)$. 
The following proposition is not very difficult and shows the analog for the case 
with boundary.   

\begin{prop}\label{1.5}{A bounded operator on $L^2(X)$ with range in 
$H^n(X)$ is an element of $\cL^{(1,\i)}(L^2 (X))$; if its range
even is contained in  $H^{n+1}(X)$ then it is trace class.
%The embedding $H^n (X) \hookrightarrow L^2 (X)$ belongs to $\cL^{1,\i}
%(L^2 (X))$, while $H^{n+1}(X) \hookrightarrow L^2(X)$ is trace class.
}
\end{prop}

\defn\label{2.0}{We shall denote by $\D^{m}$ the space of all  operator matrices 
$$A= \left(\begin{array}{cc}
P_++G & K\\
T & S
\end{array}\right) $$ in Boutet de Monvel's calculus, where 
$P$ is of order $m$, $G$ is of order $m$ and type zero,
$K$ is of order $m$, $T$ is of order $m+1$ and type zero, and $S$ is of order $m+1$. 
}\medskip

\rm 
The following lemma is now obvious from  \ref{1.7}(g):

\cor\label{2.1}{An operator $A\in \D^{-n}$ defines a bounded map
$$A:H\to H^{-n}(X)\oplus H^{1-n}(\partial X).$$
Hence ${\D}^{-n} \hookrightarrow \cL^{(1,\i)} (H) $. Similarly 
${\D}^{-n-1}  \hookrightarrow \cL^1 (H)$.
In particular, Dixmier's trace applies to elements in ${\D}^{-n}$; it vanishes on ${\D}^{-n-1}$.}\medskip
 
%\rem{\rm Obviously, no multiple of the noncommutative residue in \eqref{fgls} 
%can coincide with Dixmier's trace
%(unless $\partial X = \emptyset$): The contributions by $P$ and $S$ have
%opposite signs, while Dixmier's trace is a {\em positive} functional.}
%\medskip

\rm We will show the following theorem: 

\begin{thm}\label{2.6}{
For an operator
$$
A = \left( \begin{array}{cc} P_+ + G & K \\ T & S \end{array} \right) \in
\D^{-n}
$$
acting on $H = L^2 (X,E) \oplus L^2 (\partial X,F),\  n = {\rm dim}\, X>2$,
we have
\begin{eqnarray*}
\lefteqn{{\rm Tr}_\go (A) = \frac{1}{(2\pi)^n n} \int_X \int_S {\rm tr}_E\, p_{-n}
(x,\xi)\, \gs(\xi) dx} \\
&& {} + \frac{1}{(2\pi)^{n-1} (n-1)} \int_{\partial X} \int_{S'}
{\rm tr}_F \,s_{-n+1} (x',\xi') \,\gs'(\xi') dx' .
\end{eqnarray*}
\rm Here, $S$ is the sphere $\{|\xi|=1\}$ in $T^*\gO$ over a point
$x\in X$, $\gs$ is the form 
introduced in {\rm\ref{251.2}}, while $S'$ is the corresponding sphere 
$\{|\xi'|=1\}$ in $T^*\partial X$ over a point $x'$ in $\partial X$ and $\gs'$ 
is the corresponding $(n-2)-$form.

For ${\rm dim}\, X =2$ the second summand becomes
$$
\frac{1}{2\pi} \int_{\partial X} {\rm tr}_F (s_{-1} (x',1) + s_{-1} (x',-1))\,
dx' .
$$
}
\end{thm}\bigskip

\rm We start with the following observation.

\begin{lem}\label{2.2}{
Let $L \in {\mathcal K} (L^2 (X,E))$. Then, apart from zero, the spectrum of $L$
in $\cL (L^2 (X,E))$ including multiplicities coincides with the spectrum of
$\left( {\vec L0} {\vec00} \right)$ in $\cL
(H)$ and the spectrum of $1_X L 1_X$ in $\cL (L^2 (\gO,E))$.

{\rm Here $1_X$
denotes the characteristic function of $X$ in $\gO$, and we identify
$L^2$-functions on $X$ with their extensions by zero.}
}
\end{lem}\medskip 

We obtain a first result:

\begin{lem}\label{2.3}{
Let $\left( {\vec GT} {\vec K S} \right) \in \D^{-n}$ have
zero pseudodifferential part. Then
$$
{\rm Tr}_\go \left( \begin{array}{cc} G & K \\ T & S \end{array} \right)
 =
\frac{{\rm res}_{\partial X}\, S}{(2\pi)^{n-1} (n-1)} 
$$
with Wodzicki's residue on $\partial X$. There is no contribution from
$G,T$, or $K$.
}
\end{lem}

\Proof. \rm 
${\rm Tr}_\go$ is a trace on $\cL^{(1,\i)} (H)$, therefore ${\rm Tr}_\go ([B,C]) = 0$
whenever $B \in \D^{-n} (X)$ and $C \in \cL (H)$, 
cf.~Proposition 1.16 (d). 
In particular, we easily deduce from commutator identities  \eqref{c1} and \eqref{c2}
(with $\tau = {\rm Tr}_\go$) that Dixmier's trace vanishes on the off-diagonal entries. 
We know from (an analog of) Lemma \ref{2.2} and Theorem \ref{10.1} that
\begin{equation}\label{trs}
{\rm Tr}_\go \left(  \begin{array}{cc} 0 & 0 \\ 0 & S \end{array} 
\right) = {\rm Tr}_\go (S) = \frac{{\rm res}_{\partial X}\, S}{(2\pi)^{n-1} (n-1)} .
\end{equation}
Next let $K'$ be a potential operator of order zero. Then $K' : L^2 (\partial X,F) \to L^2 (X,E)$ and $ K':H^n(\partial X,F)\to H^n(X,E)$ are bounded by \ref{1.7}(g). 
Pick also a trace operator
$T'$ of order $-n$ and type zero.  
The operator $T' : L^2 (X,E) \to H^n (\partial X,F)$ is bounded,    and  $T'K'$  a trace class
operator on $L^2(\partial X,F)$, for
$\dim\, \partial X= n-1$.
Moreover, $B=\left( {\vec00} {\vec{K'}0} \right) \in \cL (H)$, and 
$C=\left( \vec0{T'} {\vec00}\right) \in \D^{-n}$. We conclude from \eqref{trs} that 
$$
0 = {\rm Tr}_\go \left( \left[B,C\right]
\right) = {\rm Tr}_\go \left( \begin{array}{cc} K'T' & 0 \\ 0 & -T'K'
\end{array} \right)  = {\rm Tr}_\go \left( \vec{K'T'} 0 {\vec0  0}
\right) .
$$
The span of the operators $K'T'$ of this form  is dense in
the space of singular Green operators of order $-n$ and type zero. 
As a consequence,  Dixmier's trace vanishes on
the upper left corner, since it is continuous on $\cL^{1,\i} (H)$, and
the topology on the singular Green operators is stronger than the
operator topology.
By linearity the proof is complete.
\eproof\medskip

\rm %\begin{cor}\label{2.4}{
The task of determining   Dixmier's trace for the matrix $A$ in Definition \ref{2.0}
reduces to the case where $G$, $K$, $T$, and $S$ are zero. 
By Lemma \ref{2.2} this amounts to finding 
${\rm Tr}_\go (P_+)$ in $\cL(L^2 (X,E))$ or 
${\rm Tr}_\go (1_X P1_X)$ in $\cL (L^2 (\gO,E))$. 
Since  ${\rm Tr}_\go$ is a trace, the latter equals 
$  {\rm Tr}_\go (1_X P) $. 

%Note that the fact that $P$ satisfies the transmission property is no
%longer of importance.
%}
%\end{cor}

\begin{prop}\label{2.5}{
For a  pseudodifferential operator $P$ of order $-n$, acting on $L^2(\gO,E)$, 
\begin{equation}\label{tr1xp}
{\rm Tr}_\go (1_X P) = \frac{1}{(2\pi)^n n} \int_X \int_S {\rm tr}_E p_{-n} (x,\xi)\,
\gs(\xi) dx .
\end{equation}
{\rm Here, $P$ need not have the transmission property.}
}
\end{prop}

\Proof. 
Let us first  assume  that $P$ is positive. Choose sequences of
smooth functions $\{ \gvp_k \}, \{ \psi_k \}$ with $0 \leq \gvp_k \leq
1_X \leq \psi_k \leq 1$ and $\gvp_k, \psi_k \to 1_X$ pointwise.

Since $\gvp_k, \psi_k$ are smooth, $\gvp_k P$ and $\psi_k P$ are
the pseudodifferential operators with symbols $\gvp_k p$ and $\psi_k p$,
respectively. For them, Dixmier's trace is given by Theorem \ref{10.1}. In view of
the positivity of ${\rm Tr}_\go$,  
$${\rm Tr}_\go (1_XP-\gvp_kP)={\rm Tr}_\go((1_X-\gvp_k)P)= 
{\rm Tr}_\go(\sqrt{1_X-\gvp_k}~P~\sqrt{1_X-\gvp_k})\ge0.$$
Arguing similarly for $\psi_k$, we conclude that
\begin{eqnarray*}
\frac{1}{(2\pi)^n n} \int_\gO \int_S {\rm tr}_E \gvp_k (x) p_{-n}
(x,\xi) \gs(\xi) dx&=& {\rm Tr}_\go (\gvp_k P) \leq {\rm Tr}_\go (1_X P) \leq {\rm Tr}_\go (\psi_k P)\\
&=& \frac{1}{(2\pi)^n n} \int_\gO \int_S {\rm tr}_E \psi_k (x) p_{-n} (x,\xi)\,
    \gs(\xi) dx .
\end{eqnarray*}
Lebesgue's theorem on dominated convergence then gives \eqref{tr1xp}.

Next assume that $P$ is selfadjoint. 
Let $\Delta$ be the Laplace-Beltrami operator associated with an arbitrary
Riemannian metric. Define $T = (I-\Delta)^{-n/4}$. This is a positive
selfadjoint pseudodifferential operator of order $-n/2$.
It is invertible, and $T^{-1}$ is of order $n/2$. 
The operator 
$T^{-1} P T^{-1}$  is of order zero and therefore bounded;
it also is selfadjoint. 
For large $t\in\R$, we conclude that $T^{-1} P T^{-1} + t I$ is positive. 
Hence also  $P+tT^2= T(T^{-1} P T^{-1} + t I) T $ is positive.
We can therefore write $P$ as the difference of the positive operators
$P+tT^2$ and $tT^2$. In view of the additivity of the integral
we obtain the assertion. 

Finally, in the general case we write $
P = \frac{1}{2} (P+P^*) - \frac{i}{2} (iP - iP^*) . $ 
Since $P \pm P^*$ is pseudodifferential of order $-n$, additivity  shows the assertion. 
\eproof\bigskip

The proof of Theorem \ref{2.6} is now complete. 
\bigskip

In the more traditional literature one adopts slightly different way of 
looking at boundary value problems.  
We shall now sketch how the results apply to this  situation. 

%\refstepcounter{thm}
 %{\sc \thesection.\arabic{thm} 
{Let $P:C^\i(X,E)\to C^\i(X,E)$ 
be a differential operator of  order $m>0$ and $T:C^\i(X,E)\to C^\i(Y,F)$   a trace operator. 
One usually  assumes that 
$$F=F_1\oplus \ldots \oplus F_m,$$
where the dimension of each $F_j$ might be zero (in a standard application, 
$m$ would be even and half of the $F_j$ would be zero).
Correspondingly one writes $T=(T_1, \ldots,T_m)$, asking that $T_k$ be a 
differential boundary operator of order $m_k< m$ involving at most $k-1$ normal derivatives. 
In local coordinates near the 
boundary, $T_k$ can be written in the form 
$$T_k=\sum_{j< k,~ j+|\ga|\le m_k}a_{j\ga}^{[k]}(x')D^\ga_{x'}\gg_j.$$
Here,  $\gg_j=\gg_0\partial_{x_n}^j$ 
is the operator of evaluation of the $j$-fold normal derivative at the
boundary. Of course, $T_k$ can be neglected when the dimension of $F_k$ is zero.

In order to treat the homogeneous problem
\begin{equation}\label{271}
Pu=f\ \ \ {\rm on} \ \ X;\qquad Tu=0\ \ \ {\rm on} \ \ \partial X
\end{equation} 
one studies the ``realization'' $P_T$  which is 
defined as the unbounded operator $P_T$ on $L^2(X,E)$ with domain 
\begin{equation}\label{274}
{\mathcal D}(P_T) = \{u\in H^m(X,E): Tu=0\}
\end{equation}
and acting like $P$ rather than the operator
\begin{equation}\label{273}
{\binom{P}{T}}: H^m(X,E) \to  
\begin{array}{l} L^2(X,E) \\ \oplus \\
\bigoplus_{k=1}^{m} H^{m-m_k-1/2}(\partial X,F_k).
\end{array}
\end{equation}
In both cases one is interested in the question whether the corresponding operator
has the Fredholm property. 
The link between the two approaches is given by the following proposition. }

\prop\label{2.8}{Let $H$, $H_1$ and $H_2$ be Hilbert spaces over $\C$, 
and let  $P: H\to H_1$ as well as $T: H\to H_2$ be  linear operators. 
Then the following are equivalent:
\bli
\item $\left(\begin{array}{c}P  \\ T\end{array}\right): H \to  \begin{array}{c} H_1 \\ \oplus \\H_2\end{array}$
is a Fredholm operator.
\item $P_T = P|_{{\rm ker }\, T}: {\rm ker }\, T\to H_1$  is a Fredholm operator, and 
${\rm im }\,T $ is finite codimensional in $H_2$.
\eli}

\rm In standard cases one will often have $T$ surjective so that the
second condition in (ii) is easily fulfilled. \medskip

\refstepcounter{thm}
 {\sc \thesection.\arabic{thm} Ellipticity. }\label{2.9}{
\rm The operator in 
\ref{273}
%(and therefore also \ref{2.7}(\ref{271})) 
is a Fredholm operator 
if and only if the following two  conditions are fulfilled:
\renewcommand{\labelenumi}{(\roman{enumi})}
\begin{enumerate}
\item For all $(x,\xi)\in S^*\gO|_X$ the principal symbol $p_m$ of $P$ is invertible as an
endomorphism of $E$, and
\item for all $(x',\xi')\in S^*\partial X$ the principal boundary symbol 
$$\left(\vec{p_m(x',0,\xi',D_n)}{t_{m-1/2}(x',\xi',D_n)}\right): \cS(\rp,E)
 \to  \begin{array}{c}\cS(\rp,E)\\ \oplus \\ F_{(x',\xi')} \end{array}
$$
is invertible. 
\end{enumerate}
Here $t_{m-1/2}$ is the operator-valued
principal symbol of the trace operator $T$. Locally near the boundary, this  is the vector with entries
$\sum_{j< k,~j+|\ga|= m_k}a_{j\ga}^{[k]}(x'){\xi'}^\ga \gg_j.$}

The orders of the entries in this boundary value problem are not as required in Definition 
\ref{1.6}. This, however, need not worry us:
We may replace each 
$T_k$ by $\ti T_k= \gL^{m-m_k-1/2}T_k$, 
where  $\gL = (1-\Delta_{\partial X})^{1/2}$ is order-reducing 
along $\partial X$. In view of the fact that
the powers of $\gL$ are invertible, this affects neither  
 \eqref{274} nor
the ellipticity of $\binom P T$, 
% \ref{2.7}(\ref{271}) nor the Fredholm property, 
but it 
allows us to use the space $L^2(\partial X,F)$ instead of $\oplus_{k=1}^m
H^{m-m_k-1/2}(\partial X, F_k)$ on the right hand side of \eqref{273}. 

The ellipticity implies that there is a parametrix to 
$A={\binom {P}{\ti T}}$ in Boutet  de Monvel's calculus. 
It is of the form $B=(Q_++G~~K)$; the pseudodifferential part $Q$   
is a parametrix to $P$, while $G$ is a singular Green operator of order $-m$ and 
type zero and $K$ is a potential operator of order $-m$.
Being a parametrix here means that 
$AB-I = S_1$ and $BA-I=S_2$ are regularizing operators; their 
types are $0$ and $m$, respectively. As a consequence the operator 
$S_1$ is an 
integral operator with a smooth kernel section. 

Multiplication by
${\rm diag}\,(1,\gL_1,\cdots,\gL_m)$ with $\gL_k = \gL^{m-m_k-1/2}$,
furnishes a para\-metrix for $\binom PT $;
note that  this only affects the entry  $K$ of $B$. }

The Theorem, below, follows from a construction by Grubb and Geymonat,
see e.g. \cite{Grubb86}, Section 1.4.

\thm\label{2.11}{Let $\binom{P}{T}$ be the above elliptic boundary value problem and 
$B=(Q_++G\ \ K)$ as above.
Then there is a regularizing singular Green operator $G_0$ 
of type zero, i.e.~an integral operator with smooth kernel on 
$X\times X$, such that $R= Q_++G +G_0$
has the following properties:
\def\labelenumi{(\alph{enumi})}
\begin{enumerate}
\item $R$ maps $L^2(X,E)$ to ${\mathcal D}(P_T)$ and 
\item $RP_T-I$ and $P_TR-I$ are finite rank
operators whose range consists of smooth functions. 
\end{enumerate}
Any other parametrix for $P_T$ differs from $R$ by a regularizing s.G.o..   
}\medskip

\cor\label{2.13}{\rm Suppose that  $m=n$ and $\binom{P}T$ is elliptic. Then  
Theorem \ref{2.6}   shows that Dixmier's trace 
for an arbitrary parametrix $R$ to the operator $P_T$ is given by
\begin{equation}\label{2131}
{\rm Tr}_\go R = \frac{1}{(2\pi)^nn}\int_X\int_S
{\rm tr}_E p_{n}(x,\xi)^{-1} \gs(\xi) dx.
\end{equation}
In particular, it  is the same for all parametrices  and independent 
of the  boundary condition $T$. It coincides with the noncommutative residue for $R$.
}\medskip

\rm
For certain classes of  parameter-elliptic boundary value 
problems one can deduce this result from eigenvalue asymptotics for
the operator $P_T$. Results of this type were first established by H. Weyl
\cite{Weyl11} for the 
Dirichlet problem in the plane and subsequently 
by numerous authors for many operators.
  Theorem \ref{2.6} on the other hand covers  more general 
situations, in particular, non-elliptic operators.
For example, we get the result of \ref{2.13} for all elliptic boundary value 
problems $\binom{P+G}{T}$ of order $n$ in Boutet de Monvel's 
calculus (i.e., without requiring parameter--ellipticity and with
additional singular Green terms as well as  more general trace operators);
formula (\ref{2131}) continues to hold.

\section {Heat Trace Asymptotics}
\rm
It is an obvious question whether there is an analog of the expansion \eqref{E2} 
and relation \eqref{c}. In fact one might be tempted to introduce the 
noncommutative residue for manifolds with boundary in this way. There is, however,
the following obstacle:\medskip

\refstepcounter{thm}
 {\sc \thesection.\arabic{thm} Earlier Results. }\label{ER}
\rm  Let $B=(\ti P_+ + \ti G)_{\ti T}$ be the realization of an elliptic boundary 
value problem $\binom {\ti {P_+} + \ti G}{\ti T}$ in Boutet de Monvel's 
calculus. Here, $\ti P$ and $\ti T$ are as in the previous section; $\ti G$ is a s.G.o.. 
Under suitable 
assumptions on $\ti P$, $\ti G$ and $\ti T$,  Grubb showed in her book 
\cite{Grubb86} that $e^{-tB}$ is a trace class operator on $L^2 (X,E)$ and 
studied the behavior of trace $\exp(-tB)$ as $t \to 0^+$. She obtained an 
expansion similar to \eqref{eta}. For technical reasons, however, this expansion
was finite and the remainder term
was $O(t^{1-\varepsilon})$, with $\varepsilon > 0$ determined by the so-called 
regularity number of the operator.

Given an operator $P_+ + G$ in Boutet de Monvel's calculus (for simplicity let 
$F=0$) one would have liked to choose an elliptic boundary value problem 
$\binom{\ti {P_+} + \ti G}{\ti T}$ of higher order as above and to define ${\rm 
res}\, (P_+ + G)$ by the derivative $d/du|_{u=0}$ of the coefficient of $t \ln t$ in the 
expansion of ${\rm trace }\,\exp(-t(\ti P_1 + \ti G + u(P_+ + G))_{\ti T}$ for small $|u|$, just as in  \eqref{var}. Due to the uncertainty $O(t^{1-\varepsilon})$, however, 
this does not make sense. 

The difficulties one encounters when considering trace $(P_+ + G) \exp 
(-t (\ti P_+ + \ti G)_{\ti T})$ are the same. Hence \eqref{255} can not be employed, either. 
It even does not help to restrict $\ti P$, $\ti G$,  and  $\ti T$ to suitable classes of 
operators. \medskip

Progress was possible as a result of the work of Grubb and Seeley. Already in 
\cite{GrubbSeeley93a} and \cite{GrubbSeeley93c} they applied the weakly 
parametric calculus in order to study trace expansions for operators arising 
from Atiyah-Patodi-Singer boundary value problems.

In \cite{GrubbSchrohe99} we now employ this calculus in order to 
obtain the complete expansion for 
$${\rm trace }\,(P_+ + G) \exp (-t (\ti P_+ + \ti G)_{\ti T}),$$
for a certain class of  auxiliary operators. This yields the  
 link between the noncommutative residue and 
heat trace expansions. \medskip

\refstepcounter{thm}
 {\sc \thesection.\arabic{thm} Outline. }\label{O}
We proceed in the following way.
We pick an elliptic second order pseudodiffererential operator $P_1$ on $\gO$ 
whose
principal symbol $p_2 = p_2 (x, \xi)$ has its eigenvalues in a sector in the 
right half-plane. We let $P_{1,Dir}$ be the Dirichlet realization of $P_1$, 
i.e.~the unbounded operator on $L^2(X)$ with domain
$
{\mathcal D}(P_{1,Dir}) = \{u \in H^2(X) : \gamma_0 u = 0 
$
on
$\partial X\}$.

Moreover, we choose an elliptic second order pseudodifferential operator $S_1$ 
on $\partial X$, asking again that the eigenvalues of the principal symbol lie 
in a sector of the right half-plane.

Just to give an example, we might fix arbitrary Riemannian metrics $g$ on $X$ 
and $g'$ on $\partial X$ and let $P_1 = \Delta_g {\rm Id}_E,$ $S_1 = 
\Delta_{g'} {\rm Id}_F$ be the associated Laplace-Beltrami operators.

We then let
$$
A = \left (\begin{array}{cc}P_{1,Dir}& 0 \\ 0 & S_1
\end {array} \right).
$$
{}From the analysis in \cite[Section 4.2]{Grubb86} one immediately obtains the 
following theorem. 

\thm\label{T1}{The spectrum of $A$ lies in an obtuse key hole regin
$$ W=\{\gl\in\C: |\gl|\le r~~{\rm  or}~~ |\arg \gl|\le \pi/2-\gve\},$$
where $r>0$ and $0<\gve<\pi/2$ are fixed. The
 operator $A$ generates a strongly continuous semigroup 
$e^{-tA}$ on $L^2(X,E) \oplus L^2(\partial X,F)$. For each $t>0$, the operator 
$e^{-tA}$ is given by the Dunford integral 
$$e^{tA} = \frac{i}{2\pi} \int_{\mathcal C} 
e^{-tA} (A-\lambda)^{-1} d\lambda,$$
where $\cC$ is the counter-clockwise oriented boundary of $W$.
It maps $L^2(X,E) \oplus L^2(\partial X,F)$ to 
$C^\infty(X,E) \oplus C^\infty(\partial X,F)$.
}\medskip

\rm For every operator 
$$ {\cP} = \left (\begin{array}{cc}P_+ + G & K\\ T & S \end{array} 
\right)$$ 
in Boutet de Monvel's calculus, the composition
%$\left (\begin{array}{cc}P_+ + G & K\\ T & S \end{array} \right) 
${\cP}\,e^{-tA}$    
 therefore is a trace class operator for $t>0$. 

Using identity \eqref{ETA} we obtain the following:

\cor\label{T2}
{\rm Let $\cP$ be an operator of order $\nu$ and type $d$. Assume w.l.o.g. that
$d\le \max\{\nu,0\}$. Fix $k>\frac{n+\nu}2$.
 Then
\begin{eqnarray}
 {\rm trace}\, \cP e^{-tA}
&=&\frac{ik!}{2\pi}(-t)^k \int_\cC e^{t \lambda}
{\rm trace }\, \{\cP(A-\lambda)^{-k}\}d \lambda \nonumber\\
&=& \frac{ik!}{2\pi} (-t)^k \int_\cC e^{t \lambda} {\rm trace} \{(P_+ +
G)(P_{1,Dir}-\lambda)^{-k}\} d \lambda\label{ss2}\\
&&+ \frac{ik!}{2\pi} (-t)^k \int_\cC e^{t \lambda} {\rm trace} \{S (S_1 -
\lambda)^{-k}\} d\lambda\nonumber\label{ss1}
\end{eqnarray}
We know from  Section 1 that there is an expansion
for ${\rm trace}\, S(S_1-\lambda)^{-k} $
%\sim \sum_{j=0}^\infty e_j \lambda ^
%\frac {n+{\rm ord}S-j}{2} +  \sum\limits_{l=0}^\infty (e'_l \ln \lambda
%+ e''_l) t^{-l-k},\qquad
as $|\lambda| \to \infty$
which in turn leads to an expansion of the second integral
$$
{\rm trace} \,S e^{-tS_1} \sim  \sum\limits_{j=0}^\infty \ti e_j t^\frac{j-n-\nu}{2} +   \sum\limits_{l=0}^\infty (\ti e'_l \ln t + \ti e''_l)
t^l.$$
The coefficient $\ti e_0'$ of $\ln t$ is given by
$$\ti e_0' =
\frac{-1}{2 (2\pi)^{n-1}}\,{\rm res} _{\partial X}\, S .
 $$
}\medskip

\rm It therefore remains to deal with the first integral. To this end we let
$\lambda = \mu^2$ outside $W$. We recall from \ref{2.11} that the inverse of
$P_{1,Dir}-\mu^2$ is of the form
$$
(P_1-\mu^2)^{-1}_+ + G_\mu
$$
where $\{G_\mu : \mu^2 \in \C\setminus W\}$ is a family of s.G.o.s (it can  be described
rather precisely, but we shall presently not go into the details).
We immediately deduce that
$$
(P_{1,Dir} - \mu^2)^{-k} = [(P_1 - \mu^2)^{-k}]_+ + G^{(k)}_\mu
$$
with the inverse $(P_1 - \mu^2)^{-k}$ on $L^2(\gO,E)$ and
$$
G^{(k)}_\mu = [(P_1 - \mu^2)^{-k}]_+ - [(P_1 - \mu^2)^{-1}_+]^k + R_\mu;
$$
where $R_\mu$ is a polynomial in the non-commuting operators $(P_1 - \mu^2)_+$
and $G_\mu$.

As a result, we write
\begin{eqnarray*}
(P_+ + G) (P_{1,Dir} - \mu^2)^{-k}
&=& [P(P_1-\mu^2)^{-k}]_+ - L(P, (P_1 - \mu^2)^{-k}) + P_+G^{(k)}_\mu\\
&&+ G[(P_1 - \mu^2)^{-k}]_+ + G G^{(k)}_\mu.
\end{eqnarray*}

\rm
It remains to consider the traces of these operators. The first term is easily
analyzed as a consequence of the fact that $P_1$ is an operator over
the full manifold, satisfying the assumptions in Section 1.
Recalling that relation \eqref{E2} holds
for the kernel, pointwise along the diagonal, we obtain the following assertion:

\lem\label{T3}{There is an expansion of the form {\rm \eqref{E1}} for
${\rm trace}\,P(P_1-\lambda)^{-k}$, hence an expansion of
${\rm trace}\,Pe^{-tP_1}$ of the form \eqref{E2}. The coefficient $\ti c_0'$ of $\ln t$
is given by {\rm(}recall that ${\rm ord}\, P_1=2${\rm)}
$$\ti c_0' =
\frac{-1}{2(2\pi)^{n}} \int_X \int_S {\rm tr_E}\, p_{-n} (x,\xi) \, \gs(\xi) \, dx .
 $$ }\medskip

\rm The analysis of the remaining four terms is much less trivial. It relies on a careful
study of the symbols and their compositions.  The general strategy is to reduce the
analysis to the boundary and to show that all the arizing pseudodifferential symbols
belong to the weakly parametric calculus of Grubb and Seeley, the expansion formulas
are a consequence of Theorem \ref{3.25}
We eventually show the following, provided the type is zero:

\prop\label{T4}{The traces of the operators  $L(P,(P_1 - \mu^2)^{-k})$, 
$P_+G^{(k)}_\mu$, and $G G^{(k)}_\mu$ have an expansion
$$
 \sum_{j=0}^\i a_j\mu^{n-2k+\nu-1-j}+
\sum_{l=0}^\infty (a_l'\ln \mu+a_l'')\mu^{-1-2k-l}
$$
as $|\mu|\to \infty$.
In particular, the coefficient of $\mu^{-2k}\ln \mu$ is zero. Performing the integration
in  \eqref{ss2} therefore produces no contribution to the coefficient
of $\ln t$. }

\prop\label{T5}{The trace of $G[(P_1-\mu^2)^{-k}]_+$ has an expansion
$$ \sum_{j=0}^\i b_j\mu^{n-2k+\nu-1-j}+
\sum_{l=0}^\infty (b_l'\ln \mu+b_l'')\mu^{-2k-l}.
$$
The coefficient $b'_0$ of $\mu^{-2k}\ln \mu$ is given by
$$
\frac{(-1)^k}{(2\pi)^{n-1}} \int_{\partial X}\int_{S'} {\rm tr}_E({\rm tr}\,g_{-n})(x',\xi') \gs'(\xi')dx'.
$$
}\medskip

\rm
This is precisely what we expect; note that the missing factor $2 = {\rm ord}P_1$
is due to the fact that we consider $\mu$ instead of $\gl$ and that $\ln \gl = 2\ln\mu$.
\medskip

In conclusion we obtain the following result:

\thm\label{T6}{For an operator $\cP$ of order $\nu$ and type zero in Boutet
de Monvel's calculus and for $k>(n+\nu)/2$, we have expansions
\begin{eqnarray*}
{\rm trace}\, \cP(A-\gl)^{-k}&\sim&\sum_{j=0}^\infty c_j \gl^{\frac{n+\nu-j}{2}-k} + \sum_{l=0}^\infty (c'_l\ln\gl + c''_l)\gl^{-l/2-k},\\
{\rm trace} \,\cP e^{-tA}& \sim&  \sum\limits_{j=0}^\infty \ti c_j t^\frac{j-n-\nu}{2} +   \sum\limits_{l=0}^\infty (\ti c'_l \ln t + \ti c''_l)
t^{l/2},\\
\Gamma(s) {\rm trace}\,(\cP A^{-s}) &\sim&\sum_{j=0}^\infty \frac{\ti c_j}{s+\frac{j-n-\nu}{m}} + \sum_{l=0}^\infty \left(\frac{-\ti c'_l}{(s+l/2)^{2}} +\frac{ \ti c''_l}{s+l/2}\right),
\end{eqnarray*}
for the latter assuming $A$ invertible.
The  coefficient $\ti c_0'$ of $\ln t$ satisfies the relation
\begin{multline*}
-2(2\pi)^n\, \ti c_0' =
 \int_X \int_S {\rm tr_E}\, p_{-n} (x,\xi) \, \gs(\xi) \, dx \\
+2\pi\,\int_{\partial X} \int_{S'} \left\{ {\rm tr}_E ({\rm tr}\,
g _{-n}) (x',\xi') + {\rm tr}_F s_{1-n}(x',\xi') \right\}
\sigma'(\xi') dx' \,=\, {\rm res}\, {\mathcal P}.
\end{multline*}

}
\rm As pointed out before, only the first of these expansions has to be established, the other two follow from  \eqref{AS} and \eqref{ETA}.

\section{Remarks and References to Further Work}

\rm 
In the one-dimensional case, Wodzicki's residue had been known before; it had been used
by Manin \cite{Manin78} and Adler \cite{Adler79} in their work on algebraic aspects of
the Korteweg-de Vries equation. 

In 1987 Wodzicki gave a more detailed account of the noncommutative 
residue and several related topics, among them `higher' residues cf.~\cite{Wod87}. 
A very good survey was compiled
by Kassel in 1989 for the S\'eminaire Bourbaki, \cite{Kassel89}.

Guillemin discovered the noncommutative residue independently in his so-called
`soft' proof of Weyl's formula on the asymptotic distribution of eigenvalues of  
operators \cite{Guillemin85}: For  a self-adjoint operator $P$, given as $\Op^w p$ for a classical Weyl symbol $p$,  the counting function $N_P(\lambda)$ of all eigenvalues of $P$
which are less than $\lambda$, satisfies 
$N_P(\lambda) \sim \gg\,{\rm vol}\,\{p\le \lambda\}$ as $\gl\to \infty$.
Here, $\gg$ is a universal constant and vol the symplectic volume.   

In \cite{Connes88}, Connes proved Theorem \ref{10.1} and used the coincidence
of the noncommutative residue and Dixmier's trace in  roughly the following way: 
For an algebra $\A$ and  a $p$-summable Fredholm module $(\cH,F)$ over $\A$, 
 he  introduces a `curvature' $\theta$ and an abstract Yang-Mills action $I(\theta)$,
using Dixmier's trace. 
In the case of a 4-dimensional smooth compact Riemannian 
${\rm Spin}^c$ manifold  he shows that  the classical 
Yang-Mills action $YM(A)$ for a connection $A$ can be recovered by 
$YM(A) = 16\pi^2 \inf I(\theta)$
with the infimum taken over a suitable class of connections related to $A$. 

In conformal field theory, Wodzicki's residue has been employed to construct central 
extensions of the Lie algebra of (pseudo-)differential symbols, 
cf.~Khesin and Kravchenko \cite{Khesin91}, Radul \cite{Radul91}.
It  also  has been used to derive an action for gravity in the framework of noncommutative geometry, Kalau and Walze \cite{KalauWalze93}, Kastler \cite{Kastler95}.

There are extensions in various directions. 
Guillemin \cite{Guillemin93}, \cite{Guillemin93a} defined traces on algebras of 
Fourier integral operators. Lesch \cite{Lesch97} extended the residue to symbols including logarithmic terms in their expansion. Melrose \cite{Melrose95} as well as Lesch and Pflaum \cite{LeschPflaum98} considered certain classes of 
parameter-dependent operators. 
Traces on operator algebras on manifolds with conical singularities have been studied in the paper \cite{Schrohe96b} by the author. 
While new traces arise, there also is the noncommutative residue defined on an ideal. 
Melrose and Nistor \cite{MelroseNistor96} proved an index theorem for an
algebra of cusp pseudodifferential operators on  manifolds with boundary.
They computed the  Hochschild cohomology groups and expressed the index in terms 
of various trace functionals on ideals, among them Wodzicki's. 
Their work is partly based on the computation of the homology of the algebra of pseudodifferential operators by Brylinsky and Getzler \cite{BrylinskiGetzler87}, where 
Wodzicki's residue as well as the higher analogs naturally arise.
 Kontsevich and Vishik \cite{KontsevichVishik95}, cf.~also \cite{KontsevichVishik95a},   finally introduced a trace functional
${\rm TR}$ on classical pseudodifferential operators on noninteger order.
For  a holomorphic family of classical pseudodifferential operators $\{A(z):z\in\C\}$
with ${\rm ord}\, A(z)=z$ they showed that ${\rm TR}\, A(z)$ is a meromorphic function
on $\C$ with at most simple poles in the integers and 
${\rm Res}_{z=m}\,{\rm TR}\,A(z) =- {\rm res}\,A(z).$

\end{document}